\documentclass[12pt]{amsart}
\usepackage{subfigure}
\usepackage{graphicx}
\usepackage{rotating}
\usepackage[T1]{fontenc}
\usepackage[utf8]{inputenc}
\usepackage{amssymb}
\usepackage{amsmath}    
\def\R{I\kern -0,37 em R}
\def\P{I\kern -0,37 em P}
\def\Z{I\kern -0,37 em Z}

\title{ON THE EQUIVALENCE PROBLEM FOR GEOMETRIC STRUCTURES, II}
\author{Antonio Kumpera}
\address[Antonio Kumpera]{Campinas State University\\
Campinas, SP, Brazil}
\email{antoniokumpera@hotmail.com}
\date{June 2014}


\begin{document}

\begin{abstract}
This paper is a continuation of Part I where the general setup was developed. Here we discuss the general equivalence problem for geometric structures and provide criteria for the equivalence, local and global, of transitive structures. Cartan's Flag Systems illustrate the theory as a major example and, finally, some attention though little is given to non-transitive structures with regular orbits \textit{i.e.}, intransitivity classes.
\end{abstract}

\maketitle

\section{Introduction}

In the part I, we recalled the theory of Differential Invariants as conceived by Sophus Lie and placed it in a context best suited for the study of the equivalence problem for geometric structures. It is interesting to add some more remarks proper to what we shall discuss in the sequel and, firstly, let us recall what the Finiteness Theorem for the differential invariants, proved in Part 1, brings to the equivalence problem. A germ of a structure \textit{S} being fixed, we are interested in the formal equivalence, general or restricted, of \textit{S} with the germs of neighboring structures or those eventually located on other manifolds. We are as well interested in global equivalences, these however offering several additional difficulties of quite a different nature since globalization is always a hard topological task. Let us now place ourselves back into the neighborhoods $\mathcal{W}_k$ as described in the statement of the Finiteness Theorem of the Part I. The constancy of the $k-th$ order differential invariants provides us with a necessary and sufficient condition for the $k-th$ order equivalence namely: The $k-$jet of a structure \textit{S'} is equivalent to the $k-$jet of \textit{S} if and only if the $k-th$ order differential invariants assume the same values on both jets, these two jets belonging therefore to the same orbit in $\mathcal{W}_k$. Similarly, two infinite jets of structures, elements of \textit{lim~proj} $\mathcal{W}_k$, are equivalent whenever they are so at every finite order. Since from an integer $\mu$ onwards the higher order differential invariants are all obtained by iterated formal admissible derivations of those of order $\leq \mu$, we infer that the infinitely many equivalence conditions are all consequences of those of order $\leq \mu$ and, moreover, just of a finite number of them, namely a finite fundamental system of such invariants. In the next section we shall examine the formal equivalence problem for the formally transitive structures and, in later sections, more specific types of structures will also be discussed. As for the first mentioned structures, we shall see that the formal equivalence takes place if and only if the restrictions of the differential invariants of orders $\leq\mu$ to the flows of the same order associated to these structures assume the same values and, further, the restrictions of the derived invariants to the flows of higher orders vanish. We shall see as well that when the flows of orders $\leq\mu$ associated to the germs of two structures are transverse to the trajectories - the opposite of transitivity - then a necessary and sufficient condition for the formal equivalence is the existence of a germ of local diffeomorphism $\varphi$ of \textit{P} such that the restriction of every invariant of order $\leq\mu$ to one of these jets of structure be equal to the restriction to the other jet composed with $\varphi$, the same property thereafter holding for all the invariants of higher order.

\vspace{2 mm}

\noindent
It now only remains to determine criteria under which the hypotheses $H_2$ and $H_3$, stated in Part I, section 9, are verified and it is precisely here that the specific nature of the prolongation spaces \textit{E} intervenes. As for the hypothesis $H_2$, it will be realised, for example, upon requiring the local regularity, at a given order \textit{k}, of the distributions $\Delta_k$ and $\Delta_{k,k+1}$ together with the 2-acyclicity (or even involutivity) of $\Delta_{k-1,k}$ and arguing as in Quillen's criterion (while assuming of course that 1-acyclicity already holds on a sequence of open sets $\mathcal{W}_k$, this reviving now the hypothesis $H_1$). A more elaborate analysis (see \cite{Kumpera1967}) shows the possibility of linking the $\delta-$complex associated to the kernels $\Delta_{k-1,k}$ with the corresponding $\delta-$complex associated to the kernels of the spaces $\tilde{J}_k(\mathcal{L}_{\ell})$ and consequently, in view of the prolongation spaces structure, with that constructed by means of the kernels of the spaces $J_{\ell+k}\mathcal{L}$, the later introducing however a shifting, by one unit, in the cohomology groups. This matter will be discussed in more detail in later sections. However, it is also worthwhile to mention that we can employ, in this context, the method of \textit{characters} by extending the reasoning of \cite{Kumpera1975}, \S 24, and relate the characters of the kernels of $J_{\ell+k}\mathcal{L}$ with those of $\Delta_{k-1,k}$.

\noindent
Finally concerning the hypothesis $H_3$, not much can be simplified in the specific case of prolongation spaces. Apparently, the most efficacious criteria are provided, on the one hand, by the proposition 21.5 of \cite{Kumpera1975} that involves the bracket of formal and holonomic derivations and, on the other, by the proposition 25.4 that relies on the explicit formula (25.5) which, surprisingly, can already be found in Lie's work. The hypotheses underlying the last proposition are no other than conditions of transversality, in the strict sense, of the $k-th$ order flows of local sections \textit{S} of the fibration $E~\longrightarrow~P$ with the trajectories of $\Delta_k$ such conditions being entirely antipodal to a formal transitivity hypothesis. We finally remark that it is possible to transcribe all the results stated above in the case where $(E,\pi,P,p)$ is no more than a prolongation space relative to an infinitesimal Lie pseudo-algebra $\mathcal{L}$ (\cite{Cartan1904}).

\section{Subordinate prolongation spaces and formal transitivity}

Let us begin by stating a very naive lemma, consequence of the functoriality properties of the prolongation operation, that nevertheless is at the basis of everything that follows.

\vspace{4 mm}

\newtheorem{profile}[LemmaCounter]{Lemma}
\begin{profile}
Let $(E,\pi,P,p)$ be a finite prolongation space, $\xi$ a local vector field defined on P and $\varphi$ a local diffeomorphism such that $\alpha(\xi)=\alpha(\varphi)$. Under these conditions, $(p\varphi)_*p\xi=p(\varphi_*\xi)$.
\end{profile}

\vspace{4 mm}

\noindent
Let us next assume that the prolongation space has finite order $\ell$. In order to render notations easier, we denote by $Z\mapsto Z_h$ the projection $J_k E\rightarrow J_h E$ and by $Z~\rightarrow~\overline{Z}_{k-h}$ the semi-holonomic inclusion $J_{k-h}(J_h E)$, the same notations applying as well to $\Pi_kP$. We observe that the elements of $J_{k-h}(J_h E)$ can be canonically identified with the $(k-h)-th$ order holonomic contact elements of dimension \textit{n} that are transverse to the fibration $J_hE~\longrightarrow~P$, the contact element $\overline{Z}_{k-h}$ being issued at the point $Y_h$. This being so, the groupoid $\Pi_{\ell+k+1}P$ operates to the left (or to the right if one prefers so) on

\vspace{4 mm}

\noindent
a) $J_mE$, $m\leq k+1$, by means of the prolongation spaces law, in view of the Lemma 2, section 2 in \cite{Kumpera2014} or else by the left action (1) (\textit{loc.cit.}, pg. 5, part I), and the projection $\Pi_{\ell+k+1}P~\longrightarrow~\Pi_{\ell+k}P$.

\vspace{4 mm}

\noindent
b) $J_{m-h}(J_h E)$, $m\leq k+1$, by applying twice the above Lemma.

\vspace{4 mm}

\noindent
c) $J_{\ell+k}TP$, by the standard action.

\vspace{4 mm}

\noindent
d) $J_{k-h}(TJ_h E)$, thanks to the left action $\Lambda_{k+h}$ (\textit{loc. cit.}) followed by the standard action hence, in particular, on $TJ_k E$.  

\vspace{4 mm}

\noindent
These actions are compatible with all the fibrations and it is quite evident that a) and b) are co-variant by means of the canonical inclusion. Furthermore, the preceding lemma implies that the infinitesimal prolongation operation is compatible with these actions. In fact, taking $Y\in\Pi_{\ell+k+1}P$, $j_{\ell+k}\xi(y)\in J_{\ell+k}TP$ and $Y(j_{\ell+k}\xi(y))=j_{\ell+k}\eta(y)$ then, for all $Z\in J_hE$ composable with \textit{Y}, the element transformed by \textit{Y} via the action d) of the jet $j_{k-h}(p_h\xi)(Z)$ is equal to $j_{k-h}(p_h\eta)(Y(Z))$. Moreover, we can extend the action of $\Pi_{\ell+k+1}P$ to $J_{\ell+k}TP~\times_P~J_hE$ acting separately on each factor. The previous considerations can now be summarized by the following statement (\textit{cf}. (7), \cite{Kumpera2014}, where we replace $\overline{\lambda}_k$ by $\overline{\lambda}_{k+h}$, $\ell+k$ by $\ell+k+h$ and $TJ_kE$ by $J_h(TJ_kE)$)

\newtheorem{pros}[LemmaCounter]{Lemma}
\begin{pros}
The mapping

\begin{equation*}
\lambda_{h+(k-h)}:J_{\ell+k}TP~\times_P~J_hE\longrightarrow J_{k-h}(TJ_hE)
\end{equation*}

\noindent
is a differential co-variant with respect to the action by $\Pi_{\ell+k+h}P$ or, in other words, this mapping commutes with both actions.
\end{pros}

\noindent
In particular, we derive the following conclusions: Let $Y\in\Pi_{\ell+k+1}P$, $Z\in J_{k+1}E$,  $Y(j_{\ell+k}\xi(y))=j_{\ell+k}\eta(y')$ and assume further that $p_h\xi$ is tangent of order $k-h$ to the contact element (of order $k-h+1$) $\overline{Z}_{(k+1)-h}\in J_{k+1-h}(J_hE)$. Then the vector field $p_h\eta$ is tangent of order $k-h$ to the contact element $\overline{W}_{(k+1)-h}$, $W=Y(Z)$, since \textit{Y} transforms $\overline{Z}_{(k+1)-h}$ onto $\overline{W}_{(k+1)-h}$ and $j_{k-h}p_h\xi(Z_h)$ onto $j_{k-h}p_h\eta(W_h)$.

\newtheorem{proc}[LemmaCounter]{Lemma}
\begin{proc}
The action of $\Pi_{\ell+k+1}P$ is compatible with the tangency relations between prolonged vector fields and contact elements.
\end{proc}

\vspace{4 mm}

\noindent
Let us now introduce the (eventually singular) vector bundles, with base space $J_kE$,

\begin{equation}
\textbf{R}^0_{\ell+k}=\{(Y,Z)~|~Y\in(\textbf{R}^0_{\ell+k})_Z~,~Z\in J_kE\}
\end{equation}

\noindent
where

\begin{equation*}
(\textbf{R}^0_{\ell+k})_Z=\{j_{\ell+k}\xi(y)\in J_{\ell+k}TP~|~(p_k\xi)_Z=0~,~y=\alpha(Z)\}
\end{equation*}

\noindent
and

\begin{equation}
\mathfrak{g}_{\ell+k}=ker~(\textbf{R}^0_{\ell+k}\rightarrow\textbf{R}^0_{\ell+k-1}),
\end{equation}

\vspace{4 mm}

\noindent
the choice of the notations being justified by the fact that the linear isotropy as well as the symbol of orders $\ell+k$ of a given structure \textit{S} of species \textit{E} can simply be obtained as the inverse image of these bundles via the flow $j_kS$ or, in other terms, by the relations

\begin{equation}
\textbf{R}^0_{\ell+k}(S)=(j_kS)^{-1}\textbf{R}^0_{\ell+k}\hspace{5 mm}\mathfrak{g}_{\ell+k}(S)=(j_kS)^{-1}\mathfrak{g}_{\ell+k}~.
\end{equation}

\vspace{4 mm}

\noindent
The following sequence, where the third arrow is the morphism (7) in \cite{Kumpera2014}, is of course exact:

\begin{equation}
0\longrightarrow\textbf{R}^0_{\ell+k}\longrightarrow J_{\ell+k}TP~\times_P~J_kE\xrightarrow{\overline{\lambda}_k}\Delta_k\longrightarrow 0
\end{equation}

\vspace{4 mm}

\noindent
The preceding lemmas show in particular that each term of this sequence is invariant by the corresponding action of $\Pi_{\ell+k+1}P$. Moreover, the action of $\Pi_{\ell+k+1}P$ induced on $\textbf{R}^0_{\ell+k}$ factors to $\Pi_{\ell+k}P$ since it is no other than the parts of the orders $\leq\ell+k$ that operate on the total isotropy $J^0_{\ell+k}TP$. Finally, we observe that the action of $\Pi_{\ell+k}P$ on $J^0_{\ell+k}TP~\times_P~J_kE$ leaves invariant the sub-space 

\begin{equation*}
\mathfrak{g}_{\ell+k}\subset (S^{\ell+k}T^*P~\otimes~TP)~\times_P~J_kE
\end{equation*}

\vspace{4 mm}

\noindent
and the induced action on the tensorial term factors, via $\Pi_{\ell+k}P~\longrightarrow~\Pi_1P$, since it is only the first order part that operates on the symbols.

\vspace{4 mm}

\newtheorem{proceed}[LemmaCounter]{Lemma}
\begin{proceed}
For all $Z\in J_kE$, the symbol $(\mathfrak{g}_{\ell+k})_Z$ is the $k-th$ algebraic prolongation of $(\mathfrak{g}_{\ell})_{\beta(Z)}$ and, consequently, $\mathfrak{g}_{\ell+k}$ is entirely determined by $\mathfrak{g}_{\ell}$ via the relation $\mathfrak{g}_{\ell+k}=\mathfrak{g}^{(k)}_{\ell}\otimes_E J_kE$ where $\mathfrak{g}^{(k)}_{\ell}$ is the $k-th$ algebraic prolongation of $\mathfrak{g}_{\ell}$. 
\end{proceed}

\vspace{4 mm}

\noindent
In fact, it suffices to choose a section \textit{S} of \textit{E} such that $Z=j_kS(y)$ and apply the proposition 4 in \cite{Kumpera2014}. The previous discussion also entails the

\vspace{4 mm}

\newtheorem{procedure}[LemmaCounter]{Lemma}
\begin{procedure}
Let $\Omega_0$ be an orbit of $\Pi_{\ell}P$ in E. For every k, the restriction $\mathfrak{g}_{\ell}~|~\beta^{-1}(\Omega_0)$ is a vector bundle of constant rank on the manifold $\beta^{-1}(\Omega_0)$. Further, taking any two points $Z,Z'\in\beta^{-1}(\Omega_0)$ and assuming that the jet $Y\in\Pi_{\ell}P$ transforms $\beta(Z)$ in $\beta(Z')$, then $(\mathfrak{g}_{\ell+k})_Z$ is transformed onto $(\mathfrak{g}_{\ell+k})_{Z'}$ by the tensorial extension of the $1-$jet $Y_1$, projection of Y, considered as a linear transformation $T_{\alpha(Z)}P\longrightarrow T_{\alpha(Z')}P$.
\end{procedure}

\vspace{4 mm}

\newtheorem{procedures}[CorollaryCounter]{Corollary}
\begin{procedures}
For all k, the Spencer $\delta-$cohomology complexes associated to the symbols $(\mathfrak{g}_{\ell+k})_Z$ and $(\mathfrak{g}_{\ell+k})_{Z'}$ taken at any two arbitrary points of $\beta^{-1}(\Omega_0)$ are isomorphic by means of the restriction of an isomorphism, of the total complexes, produced by the linear part of an element belonging to $\Pi_{\ell}P$ that sends $T_{\alpha(Z)}P$ into $T_{\alpha(Z')}P$. Consequently, the symbols have the same homological properties along the sub-manifold $\beta^{-1}(\Omega_0)$. In particular, when S is a section of E taking its values in $\beta^{-1}(\Omega_0)$ then the symbol $\mathfrak{g}_{\ell+k}(S)$ is a vector bundle of constant rank and each fibre has the same homological properties. 
\end{procedures}

\vspace{4 mm}

\newtheorem{process}[LemmaCounter]{Lemma}
\begin{process}
Let $\Omega_0$ be an orbit of $\Pi_{\ell}P$ in E. The for every k, the vector bundle $\mathfrak{g}^{(k)}_{\ell}~|~\Omega_0$ is of constant rank and is invariant by the action of the groupoid $\Pi_{\ell}P$ on $(S^{\ell+k}T^*P~\otimes~TP)~\times_P~\Omega_0$, the action on the first term factoring by $\Pi_{\ell}P~\longrightarrow~\Pi_1P$. Furthermore, $\mathfrak{g}_{\ell+k}~|~\beta^{-1}(\Omega_0)=\mathfrak{g}^{(k)}_{\ell}~\times_{\Omega_0}~\beta^{-1}(\Omega_0)$. Moreover, if S is a section of E taking its values in $\Omega_0$, then the groupoid $\mathcal{R}(S)$ is transitive on $\alpha(S)$ and the linear action of $\rho_1(\mathcal{R}(S))$ leaves invariant the vector bundle $\mathfrak{g}_{\ell+k}(S)$.
\end{process}

\vspace{4 mm}

\newtheorem{invariance}[LemmaCounter]{Lemma}
\begin{invariance}
Let $\textbf{F}_k$ be an orbit, in $J_kE$, of the finite or infinitesimal action (i.e., an orbit $\Omega_k$ or $\omega_k$) of $\Pi_{\ell+k}P$ and let us write $\textbf{F}_h=\rho\textbf{F}_k,\hspace{3 mm} 0\leq h\leq k$. Under these conditions:

\vspace{4 mm}

i) $\textbf{F}_h$ is a finite or infinitesimal orbit of $J_hE$ according to the nature of the corresponding $k-$orbit.

\vspace{4 mm}

ii) The action of $\Pi_{\ell+k+1}P$ on $J_{\ell+h}TP~\times_P~\textbf{F}_h$, $0\leq h\leq k$, leaves invariant the sub-bundles $\textbf{R}^0_{\ell+h}~|~\textbf{F}_h$ and $\mathfrak{g}_{\ell+k}~|~\textbf{F}_h$, the induced action then factoring to $\Pi_{\ell+k}P$, and operating transitively in the base space $\textbf{F}_h$. The isotropies $(\textbf{R}^0_{\ell+h})_Z$ and $(\textbf{R}^0_{\ell+h})_{Z'}$ as well as the symbols $(\mathfrak{g}_{\ell+h})_Z$ and $(\mathfrak{g}_{\ell+h})_{Z'}$, at any two points $Z,Z'\in\textbf{F}_h~$, are therefore isomorphic by means of the action by elements of $\Pi_{\ell+k}P$, resp. $\Pi_1P$.

\vspace{4 mm}

iii) The action of $\Pi_{\ell+k+1}P$ on $TJ_hE~|~\textbf{F}_h$, $0\leq h\leq k$, leaves invariant the sub-bundle $\Delta_h~|~\textbf{F}_h$, the fibres at any two arbitrary points being isomorphic via the above action and, moreover, $\Delta_h~|~\textbf{F}_h=T\textbf{F}_h$.

\vspace{4 mm}

iv) The restrictions $\textbf{R}^0_{\ell+h}~|~\textbf{F}_h$, $\mathfrak{g}_{\ell+k}~|~\textbf{F}_h$ and $\Delta_h~|~\textbf{F}_h=T\textbf{F}_h$ are vector bundles of constant rank.  
\end{invariance}

\vspace{4 mm}

\noindent
By all that has already been said, the proof is obvious. Let us mention however that, in view of the exactness of the sequence (4), the vector bundles $\textbf{R}^0_{\ell+h}$ and $\Delta_h$ are \textit{simultaneously} regular \textit{i.e.}, when one is regular then the other is also necessarily regular. We next state (without proof!) a rather long lemma that outlines all the basic techniques to be used from here onwards.

\vspace{4 mm}

\newtheorem{invariant}[LemmaCounter]{Lemma}
\begin{invariant}
The hypotheses being those of the preceding lemma, for each solution S of $\textbf{F}_k$ namely, a section S of the prolongation space E such that $j_kS$ takes its values in $\textbf{F}_k$, the equations $\mathcal{R}_{\ell+h}(S)$ and $\textbf{R}_{\ell+h}(S)$, $0\leq h\leq k$, are transitive in $\alpha(S)$ and the bundles $\textbf{R}_{\ell+h}(S)$,  $\textbf{R}^0_{\ell+h}(S)$ and $\mathfrak{g}_{\ell+h}(S)$ are all of constant rank. Further, $\mathcal{R}_{\ell+h}(S)$ is a closed and locally trivial Lie (Differentiable) sub-groupoid, its differentiable structure being regularly embedded in  $\Pi_{\ell+h}\alpha(S)$. The non-linear isotropy $\mathcal{R}^0_{\ell+h}(S)$ is a closed and locally trivial sub-bundle of $\Pi^0_{\ell+h}\alpha(S)$, each fibre being a closed Lie sub-group of the total group. The projection $\rho_{\ell+h-1}\mathcal{R}_{\ell+h}(S)$ is a locally trivial Lie sub-groupoid of $\Pi_{\ell+h-1}\alpha(S)$ and the non-linear symbol

\begin{equation*}
\mathfrak{g}_{\ell+h}(S)=ker(\mathcal{R}_{\ell+h}(S)\longrightarrow\rho_{\ell+h-1}\mathcal{R}_{\ell+h}(S)) 
\end{equation*}

\vspace{4 mm}

\noindent
is a locally trivial affine sub-bundle (unless when $\ell=1$ and $h=0$, in which case the fibre is a linear group) of the total symbol above $\rho_{\ell_h-1}\mathcal{R}_{\ell+h}(S)$. Finally, the groupoid $\mathcal{R}_{\ell+h}(S)$ leaves invariant the sub-bundles $\textbf{R}^0_{\ell+h}(S)$ and $\mathfrak{g}_{\ell+h}(S)$ via the standard action on $J^0_{\ell+h}TP$ and $S^{\ell+h}T^*P~\otimes~TP$ respectively. Inasmuch, $\mathcal{R}_{\ell+h}(S)$ leaves invariant the isotropy  $\mathcal{R}^0_{\ell+h}(S)$ via the action on $\Pi^0_{\ell+h}P$ defined by conjugation, as well as the non-linear symbol  $\mathfrak{g}_{\ell+h}(S)$ via the standard action on the total symbol defined by the translations to the left or to the right. We deduce that any two fibres of $\mathcal{R}^0_{\ell+h}(S)$ and $\mathfrak{g}_{\ell+h}(S)~|~Id=\alpha(S)$ are always isomorphic as Lie groups and that those of $\mathfrak{g}_{\ell+h}(S)$ are also isomorphic as homogeneous spaces.
\end{invariant}

\vspace{4 mm}

\newtheorem{invariants}[LemmaCounter]{Lemma}
\begin{invariants}
The data as well as the hypotheses being those of the preceding Lemma, we also have the following equalities, for $0\leq h\leq k-1$, $\textbf{R}^0_{\ell+h+1}(S)=p\textbf{R}^0_{\ell+h}(S)$ and $\mathcal{R}_{\ell+h+1}(S)=p\mathcal{R}_{\ell+h}(S)$ where $p$ is the standard prolongation operator for partial differential equations.
\end{invariants}

\vspace{4 mm}

\noindent
We finally introduce the (eventually singular) vector bundle with basis $J_{k+1}E$

\begin{equation}
\textbf{R}_{\ell+k}=\{(Y,Z)~|~Y\in(\textbf{R}_{\ell+k})_Z~,~Z\in J_{k+1}E\}
\end{equation}

\vspace{2 mm}

\noindent
where

\begin{equation*}
(\textbf{R}_{\ell+k})_Z=\{j_{\ell+k}\xi(z)\in J_{\ell+k}TP~|~(p_k\xi)_Z=\overline{Z}_{(k+1)-1}\}. 
\end{equation*}

\vspace{4 mm}

\noindent
Observe that $\overline{Z}_{(k+1)-1}\in J_1(J_kE)$ is a first order contact element namely, a transverse vector sub-space of $T_{Z_k}J_kE$. The choice of these notations is justified by the relation

\begin{equation}
\textbf{R}_{\ell+k}(S)=(j_{k+1}S)^{-1}\textbf{R}_{\ell+k}.
\end{equation}

\vspace{4 mm}

\newtheorem{info}[LemmaCounter]{Lemma}
\begin{info}
The fibre bundle $\textbf{R}_{\ell+k}$ is invariant under the action of $\Pi_{\ell+k+1}P$ on $J_{\ell+k}TP~\times_P~J_{k+1}E$. If $\textbf{F}_{k+1}$ is a finite or infinitesimal orbit in $J_{k+1}E$, then the restriction $\textbf{R}_{\ell+k}~|~\textbf{F}_{k+1}$ is of constant rank and any two fibres are isomorphic by the action of some element in $\Pi_{\ell+k+1}P$.
\end{info}

\vspace{4 mm}

\noindent
Let us recall (\textit{cf.} the Theorem 2 in \cite{Kumpera2014}) that each trajectory $~\Omega_k~$ 
is a locally trivial sub-bundle of $~J_kE~$ isomorphic to the quotient space $(\Pi_{\ell+k}P)_y~/~H(Y)$. Inasmuch, an infinitesimal trajectory $\omega_k$ is a locally trivial sub-bundle since it is isomorphic to the quotient of the connected component of the unit \textit{y} in $(\Pi_{\ell+k}P)_y$ by the isotropy (we shall be careful to assume \textit{P} connected otherwise $\omega_k$ is just a bundle whose basis is the connected component of \textit{P} containing the point \textit{y}). We can thus consider any trajectory $\textbf{F}_k$, finite or infinitesimal, as a non-linear partial differential equation (a non-linear differential system) of order \textit{k} on the fibration (and prolongation space) $\pi:E~\longrightarrow~P$ that we shall call \textit{a fundamental equation of order k}.

\vspace{4 mm}

\newtheorem{germ}[LemmaCounter]{Lemma}
\begin{germ}
Let $\textbf{F}_{k+1}$ be a fundamental equation, $\textbf{F}_k$ its projection at $k-th$ order and assume that $\textbf{F}_{k+1}\cap p\textbf{F}_k\neq\phi$. Under these conditions, $\textbf{F}_{k+1}\subset p\textbf{F}_k$ and therefore $p\textbf{F}_k$ is a union of trajectories.  
\end{germ}

\vspace{4 mm}

\newtheorem{inform}[LemmaCounter]{Lemma}
\begin{inform}
A fundamental equation $\textbf{F}_k$ is $1-$integrable namely, the property $\rho_k(p\textbf{F}_k)=\textbf{F}_k$ holds, if and only if $p\textbf{F}_k\neq\phi$ (for every  fibre). This condition being verified, the prolongation $p\textbf{F}_k$ is a locally trivial affine sub-bundle of $J_{k+1}E\longrightarrow\textbf{F}_k$ where E is assumed to be a finite prolongation space.
\end{inform}

\vspace{4 mm}

Let us denote by $p_h\textbf{F}_k$ the prolongation of order \textit{k} of $\textbf{F}_k$. Arguing  as previously, we can also show the

\vspace{4 mm}

\newtheorem{informer}[CorollaryCounter]{Corollary}
\begin{informer}
A fundamental equation $\textbf{F}_k$ verifies the extended equality $\rho_k(p_h\textbf{F}_k)=\textbf{F}_k$ if and only if $p_h\textbf{F}_k\neq 0$ and, whenever this condition is verified, the following properties also hold:

\vspace{4 mm}

i) For every $\mu\leq h$,  $p_{\mu}\textbf{F}_k$ is a locally trivial (and non void) sub-bundle of $J_{k+\mu}E~\longrightarrow~\textbf{F}_k$ saturated by the orbits of order $k+\mu$, finite or infinitesimal according to the nature of $\textbf{F}_k$.

\vspace{4 mm}

ii) $p_{\mu}(p_{\eta}\textbf{F}_k)=p_{\mu+\eta}\textbf{F}_k$, $\mu+\eta\leq h$.
\end{informer}

\vspace{4 mm}

\noindent
We observe, however, that the property $p_h\textbf{F}_k\neq 0$ does not imply, in general, the $h-$integrability of $\textbf{F}_k$ namely, the property

\begin{equation*}
\rho_{\mu}(p_{\mu+1}\textbf{F}_k)=p_{\mu}\textbf{F}_k,~\mu\le h. 
\end{equation*}

\vspace{4 mm}

\noindent
The above property only implies the $1-$integrability.

\vspace{4 mm}

Let us next examine the tangent symbol of a fundamental equation. On the one hand, it is clear that this symbol at the point $Z\in J_kE$ is simply the contact element

\begin{equation*}
(\Delta_{k-1,~k})_Z=ker~[(\Delta_k)_Z\longrightarrow(\Delta_{k-1})_{\rho_{k-1}Z}],
\end{equation*}

\begin{equation*}
(\Delta_{k-1,~k})_Z\subset [(S^kT^*P~\otimes~VE)~\times_E~J_kE]_Z=[S^kT^*_{\alpha Z}P~\otimes~V_{\beta Z}E]~\times~\{Z\},
\end{equation*}

\vspace{4 mm}

\noindent
and, on the other, we see immediately that the diagram below is commutative and exact (\textit{cf.}, the second diagram at the outset of section 3 in \cite{Kumpera2014}). Curiously enough, this diagram cannot be completed everywhere by surjectivities, as is so common in commutative diagrams, since structures do not behave always as nicely as one would like to and the problem resides in the lack of \textit{h-integrability} or, stated equivalently, in the difficulty of approximating what we really look for namely, \textit{the integrability}.

\vspace{20 mm}

\newtheorem{informing}[LemmaCounter]{Lemma}
\begin{informing}
The first vertical sequence in $(7)$ is surjective at the end if and only if the first horizontal sequence verifies the same property.
\end{informing}

\vspace{4 mm}

\begin{equation}
0\hspace{51 mm} 0\hspace{49 mm} 0
\end{equation}

\begin{equation*}
\hspace{6 mm}\downarrow\hspace{51 mm}\downarrow\hspace{49 mm}\downarrow
\end{equation*}

\begin{equation*}
0\longrightarrow\mathfrak{g}_{\ell+k}\hspace{8 mm}\longrightarrow\hspace{8 mm} (S^{\ell+k}T^*P~\otimes~TP)~\times_P~J_kE\hspace{8 mm}\xrightarrow{\ell_k}\hspace{8 mm} \Delta_{k-1,~k}\hspace{10 mm}
\end{equation*}
\hspace{114 mm}$\downarrow$
\begin{equation*}
\hspace{6 mm}\downarrow\hspace{51 mm}\downarrow\hspace{45 mm}|\hspace{1 mm}\downarrow
\end{equation*}
\hspace{114 mm}$\downarrow$
\begin{equation*}
0\longrightarrow\textbf{R}^0_{\ell+k}\hspace{12 mm}\longrightarrow\hspace{12 mm} J_{\ell+k}TP~\times_P~J_kE\hspace{13 mm}\xrightarrow{\lambda_k}\hspace{11 mm}|\hspace{1 mm}\Delta_k\hspace{5 mm}\longrightarrow\hspace{5 mm} 0
\end{equation*}
\hspace{114 mm}$\downarrow$
\begin{equation*}
\hspace{7 mm}\downarrow\rho_{\ell+k-1}\hspace{38 mm}\downarrow\hspace{44 mm}|\hspace{1 mm}\downarrow
\end{equation*}
\hspace{20 mm}  $\downarrow^{------\leftarrow--------\leftarrow--------\leftarrow--------\leftarrow----}$
\begin{equation*}
0\rightarrow\textbf{R}^0_{\ell+k-1}~\times_P~J_kE\hspace{3 mm}\longrightarrow\hspace{3 mm} J_{\ell+k-1}TP~\times_P~J_kE\hspace{3 mm}\xrightarrow{\lambda_{k-1}}\hspace{3 mm}\Delta_{k-1}~\times_{J_{k-1}E}~J_kE\hspace{3 mm}
\longrightarrow\hspace{3 mm} 0
\end{equation*}

\begin{equation*}
\hspace{62 mm}\downarrow\hspace{48 mm}\downarrow
\end{equation*}

\begin{equation*}
\hspace{62 mm} 0\hspace{48 mm} 0
\end{equation*}

\vspace{5 mm}
\noindent
In fact, the above \textit{snake arrow} tells us that the following sequence is exact:

\begin{equation*}
0\rightarrow (S^{\ell+k}T^*P~\otimes~TP)~\times_P~J_kE\xrightarrow{\ell_k}\Delta_{k-1,~k}\rightarrow coker~\rho_{\ell+k-1}\rightarrow 0
\end{equation*}

\vspace{4 mm}

\noindent
The Lemma 1 can now be extended to the following assertion:

\vspace{4 mm}

\newtheorem{deforming}[LemmaCounter]{Lemma}
\begin{deforming}
Let $(E,\pi,P,p)$ be a finite prolongation space of order $\ell,~\zeta$ a local vector field defined on E, $\varphi$ a local diffeomorphism of P and we assume that $\alpha(p\varphi)=\alpha(\zeta)$. Then, $[p_k(p\varphi)]_*p_k\zeta=p_k[(p\varphi)_*\zeta]$.
\end{deforming}

\vspace{4 mm}

\noindent
We infer the following corollaries by taking also into account the Theorem 13.1 as well as the Proposition 14.2 in \cite{Kumpera1975}.

\vspace{4 mm}

\newtheorem{deform}[LemmaCounter]{Lemma}
\begin{deform}
The groupoid $\Pi_{\ell+k+1}P$ operates differentiably on the spaces $\tilde{J}_kTE$, $TJ_kE$ and $(S^kT^*P~\otimes~VE)~\times_E~J_kE$. The Lie sub-fibration $\tilde{J}_kVE$ as well as the vector sub-bundle $VJ_kE$ are invariant by this action and the following exact sequences

\begin{equation*}
0\longrightarrow J^0_kTP~\times_P~J_kE\xrightarrow{\Sigma_k}\tilde{J}_kTE\xrightarrow{p_k}TJ_k\longrightarrow 0
\end{equation*}

\vspace{2 mm}

\noindent
and

\vspace{2 mm}

\begin{equation*}
0\longrightarrow(S^kT^*P~\otimes~VE)~\times_E~J_kE\xrightarrow{\epsilon_k} TJ_kE\xrightarrow{T\rho_{k-1,k}}
\end{equation*}
\begin{equation*}
TJ_{k-1}E~\times_{J_{k-1}E}~J_kE\longrightarrow 0
\end{equation*}

\vspace{4 mm}

\noindent
are co-variant. Furthermore, the action on the term $S^kT^*P~\otimes~VE$ factors via $\rho_{\ell+1}:\Pi_{\ell+k+1}P~\longrightarrow~\Pi_{\ell+1}P$ and that on the term $J^0_kTP$ via $\rho_k:\Pi_{\ell+k+1}P~\longrightarrow~\Pi_kP$.
\end{deform}

\vspace{4 mm}

\newtheorem{deformer}[LemmaCounter]{Lemma}
\begin{deformer}
Let $\textbf{F}_k$ be a fundamental equation and $\textbf{F}_{k-1}$ its projection. By means of the canonical identification,

\begin{equation*}
\Delta_{k-1,~k}|\textbf{F}_k=ker~[T\textbf{F}_k\longrightarrow T\textbf{F}_{k-1}] 
\end{equation*}

\vspace{4 mm}

\noindent
is a locally trivial vector sub-bundle of $S^kT^*P~\otimes~VE)~\times_E~\textbf{F}_k$ invariant under the action of $\Pi_{\ell+k}P$. When $k\geq 1$ and when $Z,Z'\in\textbf{F}_k$ and $Y\in\Pi_{\ell+k}P$ are such that $Y(Z)=Z'$, then Y transforms $(\Delta_{k-1,~k})_Z$ in $(\Delta_{k-1,~k})_{Z'}$ by means of the tensorial extension of the $(\ell+1)-$jet $Y_{\ell+1}$, projection of \textit{Y}, considered as a linear map $T_{\beta(Z)}E~\longrightarrow~T_{\beta(Z')}E$. In the case where $k=0$,

\begin{equation*}
\Delta_{-1,~0}|\textbf{F}_0=ker~[T\pi:T\textbf{F}_0\longrightarrow TP]\subset VE|\textbf{F}_0,~Y\in\Pi_{\ell}P,
\end{equation*}

\vspace{4 mm}

\noindent
and we can take any jet $Y_{\ell+1}\in\Pi_{\ell+1}P$ projecting upon \textit{Y}. When $\textbf{F}_k$ is an infinitesimal trajectory, we shall take care to restrict the action just to the open subset of $\Pi_{\ell+k}P$ that leaves invariant $\textbf{F}_k$.
\end{deformer}

\vspace{5 mm}

We next consider a pair $\textbf{F}_{k+1}$ and $\textbf{F}_k$ of fundamental equations, $k\geq 0$, and verifying $\rho_k(\textbf{F}_{k+1})=\textbf{F}_k$. We shall study the geometrical properties of the fibration $\textbf{F}_{k+1}~\longrightarrow~\textbf{F}_k$ inasmuch as a sub-fibration of the affine bundle $J_{k+1}E\longrightarrow\textbf{F}_k$ above the base space $\textbf{F}_k$ (more precisely, it would be convenient to replace the term $J_{k+1}E)$ by $\rho^{-1}_k(\textbf{F}_k)\subset J_{k+1}E$ where $\rho_k:J_{k+1}E~\longrightarrow~J_kE$). Since $\textbf{F}_{k+1}$ is an orbit, this fibration is obviously a locally trivial sub-bundle, any two fibres being isomorphic by the action of the elements belonging to $\Pi_{\ell+k+1}P$. We next recall (\cite{Kumpera1975}, \S 19) that $J_{k+1}E~\longrightarrow~J_kE$ is an affine bundle, the underlying vector bundle being equal to $(S^{k+1}T^*P~\otimes~VE)~\times_E~J_kE$, the affine action being defined as follows: The point $Z\in J_kE$ being fixed, a vector $v\in S^{k+1}T^*_yP~\otimes~V_zE$, $y=\alpha(Z),~z=\beta(Z)$, defines by means of the canonical identification (\cite{Kumpera1975}, Proposition 14.2) a vector field $\nu(v)$ along the fiber $A_Z$ of $J_{k+1}E$ above the point \textit{Z}. This vector field admits a global 1-parameter group $(\varphi_t)$ of affine  transformations and the affine action of \textit{v} on $W\in A_Z$ is defined by $W+v=\varphi_1(W)$. We thus see that the set $\{\nu(v)~|~v\in S^{k+1}T^*_yP~\otimes~V_zE\}$ is the abelian Lie algebra of all the infinitesimal translations of the affine space $A_Z$. If we examine again the diagram (14.1) in \cite{Kumpera1975}, we shall observe that the pair $(v,W)$ identifies to an element in $\tilde{J}_{k+1}VE$ that projects onto $(0,Z)$ by the morphism $\tilde{\rho}_{k,~k+1}$ (see the first line of the diagram (14.1)). This means however, in setting $W=j_{k+1}\sigma(y)$ and $Z=\rho_k(W)$, that $(v,W)$ identifies to a jet $j_{k+1}(\zeta\circ\sigma)(y)$ where $\zeta$ is a vertical vector field on \textit{E} null to order \textit{k} (i.e., it vanishes up to order \textit{k}) along the image of the section $\sigma$ namely, $j_k(\zeta\circ\sigma)(y)=0~$.  Equivalently, $\zeta$ is a vertical vector field, on \textit{E, tangent} to order \textit{k} along the image of $\sigma$ and at the point $\sigma(y)$. This being so, we see readily, and according to the diagram $(14.1)$, that $\nu(v)_W=(p_{k+1}\zeta)_W~$.

\vspace{2 mm}

Let us now return to the fibration $\rho_k:\textbf{F}_{k+1}~\longrightarrow~\textbf{F}_k$. Since $\textbf{F}_{k+1}$ is an orbit under the action of $\Pi_{\ell+k+1}P$ on $J_{k+1}E$, the fibre $(\textbf{F}_{k+1})_W,~W\in\textbf{F}_k$, is the orbit of an arbitrary point $W'\in(\textbf{F}_{k+1})_W$ under the action of the sub-group $K(W)$ of $(\Pi^0_{\ell+k+1}P)_y$, inverse image of the isotropy $H(Z)$ by the projection $\rho_{k,~k+1}~$. Consequently, the fibre $(\textbf{F}_{k+1})_W$ is isomorphic to the homogeneous space $K(W)/H(Z)$ that, in general, is not connected. However, we can easily determine its connected components since they are simply the orbits of the connected component of the unity in $K(W)$ or, equivalently, the orbits of the infinitesimal action of the Lie algebra $k(W)$ of the group $K(W)$. We shall therefore study the algebra $q(W)$ of the vector fields of $A_Z$, images by the infinitesimal action of $k(W)$. We first remark that $T_Z(\textbf{F}_{k+1})_W=(\Delta_{k,~k+1})_Z$ and that consequently, every vector $w\in  T_Z(\textbf{F}_{k+1})_W$ is obtained as follows: We take a local vector field $\xi$ of \textit{P} defined in a neighborhood of \textit{y} and consider its prolongation $p_{k+1}\xi=p_{k+1}(p\xi)$. Under these conditions, $(p_{k+1}\xi)_W\in (\Delta_{k,~k+1})_W$ if and only if $T\rho_k(p_{k+1}\xi)_Z=(p_k\xi)_Z=0,~Z=\rho_k(W)$ which means precisely that $j_{\ell+k+1}\xi(y)\in k(W)$. Moreover, the condition $(p_k\xi)_Z=0$ also means, in writing $Z=j_k\sigma(y),~z=\sigma(y)$, that (\textit{c.f.}, \cite{Kumpera1975}, p.317, the remark just following the Corollary)

\vspace{4 mm}

a) $p\xi(z)=0$ and

\vspace{2 mm}

b) $p\xi$ is tangent to order \textit{k} along $im~\sigma$ at the point \textit{z}.

\vspace{4 mm}

\noindent
Since we have nullity\footnote{\textbf{Nullité.} \textit{Vice qui ôte à un acte toute sa valeur}. Larousse} at the point \textit{z}, we infer that the $k-th$ order tangency condition along the section $\sigma$ and at the point \textbf{z} only depends on the $k-th$ order jet of $\sigma$ at the point \textit{y} or, in other terms, $p\xi$ will also be tangent up to order \textit{k} along another section $\tau$ and at the point \textbf{z} as soon as $j_k\tau (y)=j_k\sigma (y)$. The vector field $\xi$ determines, by $k-th$ order prolongation, a vector field $p_{k+1}\xi$ tangent to $A_Z$ and we shall write

\begin{equation*}
\tilde{\nu}(\xi)=p_{k+1}\xi|A_Z=p_{k+1}(p\xi)|A_Z.
\end{equation*}

\vspace{4 mm}

\noindent
Since the prolongation operation preserves the brackets \textit{i.e.},

\begin{equation*}
[p_{k+1}(p\xi)~,~p_{k+1}\zeta]=p_{k+1}([p\xi,\zeta]),
\end{equation*}

\vspace{4 mm}

\noindent
where $\zeta$ is an arbitrary local vector field on \textit{E}, then

\vspace{4 mm}

1) $[p_{k+1}\xi~,~p_{k+1}\zeta]|A_Z=[\tilde{\nu}(\xi)~,~\nu(v)]=p_{k+1}([p\xi,\zeta])|A_Z$,

\vspace{4 mm}

2) $[p\xi,\zeta]$ is a vertical vector field on \textit{E},

\vspace{4 mm}

3) $[p\xi,\zeta]$ is null up to order \textit{k} along the section $\sigma$ at the point $\sigma(y)$ since $\zeta$ verifies the same property and $p\xi$ is tangent to order \textit{k} along $\sigma$ and vanishes at the point $\sigma(y)$.

\vspace{4 mm}

\noindent
We then infer that $[\tilde{\nu}(\xi)~,~\nu(v)]=\nu(v')$ where $v'\in S^{k+1}T^*_yP~\otimes~V_zE$ is the element that, for all $W\in A_Z$, determines the identification

\begin{equation*}
(v',W)\equiv (j_{k+1}([p\xi,\zeta]\circ\tau)(z)\in J_{k+1}VE~,~W=j_{k+1}\tau(z)).
\end{equation*}

\vspace{4 mm}

\noindent
The above relation shows that any vector field $\tilde{\nu}(\xi)$ belongs to the normalizer, in $W(A_Z)$, of the algebra of infinitesimal transformations of the affine space $A_Z$.Consequently, $q(W)$ is a finite dimensional Lie algebra, sub-algebra of the affine infinitesimal transformations of $A_Z$. Since the above argumentation only puts in evidence the infinitesimal transformations, we see readily that the preceding results remain valid when $\textbf{F}_{k+1}$ and $\textbf{F}_k$ are replaced just by infinitesimal orbits.

\vspace{4 mm}

\newtheorem{transformer}[TheoremCounter]{Theorem}
\begin{transformer}
Let $(\textbf{F}_{k+1}~,~\textbf{F}_k)$ be a pair of fundamental equations verifying $\rho_k\textbf{F}_{k+1}=\textbf{F}_k$. The connected components of the fibre $(\textbf{F}_{k+1})_Z~,~Z\in\textbf{F}_k$, are the orbits of the finite action, on $A_Z$, of the connected affine Lie group whose Lie algebra is equal to $q(Z)$. If $Y\in\Pi_{\ell+k+1}P$ and if $\rho_kY(Z)=W$, then $Y[(\textbf{F}_{k+1})_Z]=(\textbf{F}_{k+1})_W$ and the affine transformation $(J_{k+1}E)_Z~\longrightarrow~(J_{k+1}E)_W$ , induced by Y, transforms $q(Z)$ in $q(W)$ conjugating the corresponding affine groups.
\end{transformer}

\vspace{4 mm}

\noindent
Let us next assume that the first line of the diagram (7) is surjective at the end at a point $Y\in(\textbf{F}_{k+1})_Z~$. Since $\textbf{F}_{k+1}$ is an orbit, the Lemma 1 implies that the same property will hold in any other point $Y'\in(\textbf{F}_{k+1})_Z$. Furthermore, the surjectivity of $(\ell_{k+1})_Y$ shows that any vector in $T_Y(\textbf{F}_{k+1})_Z=(\Delta_{k,~k+1})_Y$ can be obtained by prolongation of a vector field $\xi$ satisfying the property $j_{\ell+k+1}\xi(y)\in S^{\ell+k+1}T^*_yP~\otimes~T_yP$. On account of the two commutative diagrams in the beginning of the section 3 in Part I (\cite{Kumpera2014}), we finally infer that $p_{k+1}\xi|A_Z$ \textit{is an infinitesimal translation}.

\vspace{4 mm}

\newtheorem{transformers}[CorollaryCounter]{Corollary}
\begin{transformers}
When $\ell_{k+1}$ is surjective at a point $W\in\textbf{F}_{k+1}$ then it is also surjective at any other point of $\textbf{F}_{k+1}$ and the Lie algebra $q(Z)~,~Z\in\textbf{F}_k$, contains a sub-algebra $t(Z)$ of infinitesimal translations whose orbit issued from a point of $(\textbf{F}_{k+1})_Z$ is equal to that of $q(Z)$. Each connected component of $(\textbf{F}_{k+1})_Z$ is an affine sub-space of $A_Z$.
\end{transformers}

\vspace{4 mm}

\newtheorem{transforming}[CorollaryCounter]{Corollary}
\begin{transforming}
When $\ell_{k+1}$ is surjective at a point in $\textbf{F}_{k+1}$ and if further the fibre $(\textbf{F}_{k+1})_Z$ is connected, then $\rho_k:\textbf{F}_{k+1}\longrightarrow\textbf{F}_k$ is a locally trivial affine sub-bundle of $J_{k+1}E\longrightarrow\textbf{F}_k$, any two fibres being linearly isomorphic under the action of $\Pi_{\ell+k+1}P$.
\end{transforming}

\vspace{4 mm}

\newtheorem{fake}[LemmaCounter]{Lemma}
\begin{fake}
Let $(\textbf{F}_k)_{k\geq\mu}$ be a family of fundamental equations verifying the following properties: $\rho_k\textbf{F}_{k+1}=\textbf{F}_k$ and $\textbf{F}_{k+1}\cap p\textbf{F}_k\neq\phi$. Under these conditions, $\textbf{F}_{k+1}\subset p\textbf{F}_k$ for all k and there exists an integer $\mu_0$ such that $\textbf{F}_{k+1}$ is an open sub-bundle of $p\textbf{F}_k~,~k\geq\mu_0$. The integer $\mu_0$ is the order of stability (1-acyclicity) of the Spencer $\delta-$complex associated to the tangent symbols of the equations $\textbf{F}_k$ (and their prolongations).
\end{fake}

\vspace{4 mm}

\noindent
The proof is always essentially the same and so will be omitted.

\vspace{4 mm}

We now recall that a fundamental equation does not have necessarily any solution at all. Fortunately or not, there are \textit{non-analytic} though \textit{involutive} differential systems that do not possess any solution! As for the fundamental equations, if ever they possess a germ of a solution passing by one of their points then the same will hold for all the other points and the equation will be \textit{completely (or everywhere) integrable}. It is worthwhile to state the following Corollary, consequence of the preceding Lemma, though containing some repetitions.

\vspace{4 mm}

\newtheorem{fakes}[CorollaryCounter]{Corollary}
\begin{fakes}
If $(\textbf{F}_k)_{k\geq\mu}$ is a family of integrable fundamental equations and if further the following property holds: $\rho_k\textbf{F}_{k+1}=\textbf{F}_k$, then $\textbf{F}_{k+1}\subset p\textbf{F}_k$ for all k and $\textbf{F}_{k+1}$ is an open sub-bundle of $p\textbf{F}_k$ for $k\geq\mu_0~$. We can even assume that $\mu_0=0$ since the projection of an integrable fundamental equation is also integrable.  
\end{fakes}

\vspace{4 mm}

\noindent
Let us now return to the study of \textit{symbols}. We know that $\mathfrak{g}_{\ell+k+1}$ is always the algebraic prolongation of $\mathfrak{g}_{\ell+k}$ but this is not the case, in general, for the kernels $\Delta_{k-1,k}$. Nevertheless,

\vspace{4 mm}

\newtheorem{faking}[LemmaCounter]{Lemma}
\begin{faking}
Let $(Z_k)$ be an element of $J_{\infty}E$ and let us assume, for $k\geq\mu$, that the morphism $\ell_k$ in the first horizontal sequence of $(7)$ is surjective at the point $Z_k$. Under these conditions,

\vspace{4 mm}

1) $\Delta_{k,k+1}\subset p\Delta_{k-1,k}$ for $k\geq\mu$ and

\vspace{4 mm}

2) The Spencer $\delta-$complex corresponding to the family $(\Delta_{k-1,k})_{k\geq\mu}$ is $p-$acyclic at order $\mu_p$ (i.e., for $k\geq\mu_p$) if and only if the complex constructed with the aid of the family $(\mathfrak{g}_{\ell+k})_{k\geq 0}$ is $(p+1)-$acyclic at order $\ell+\mu_p-1$
\end{faking}

\vspace{4 mm}

\noindent
As for the proof, it will suffice to recall (\cite{Kumpera2014}, sect.3) that the family $(\ell_k)$ is a natural transformation of the corresponding $\delta-$cohomology complexes compatible with the algebraic prolongation operations performed on the principal parts and proceed with the \textit{diagram chasing} below by confronting the term $\wedge^{q+1}T^*P\otimes\mathfrak{g}_{\ell+k-1}$ with $\wedge^q T^*P\otimes \Delta_{k-1,~k}$. As for the notations, we suggest to have a glance at \cite{Kumpera1972}, p.83.

\vspace{4 mm}

\noindent
An entirely analogous argument, where we replace $\Delta_{k-1,~k}$ by $h_k$, shows as well the following lemma that transcribes what we can expect when the surjectivity only occurs at the level $\mu$.

\vspace{4 mm}

\begin{sidewaysfigure}[htp!]
\vspace{12cm}
\begin{equation*}
\hspace{2 mm}0\rightarrow\wedge^{q-1}T^*P\otimes\mathfrak{g}_{\ell+k+1}\rightarrow\wedge^{q-1}T^*P\otimes(S^{\ell+k+1}T^*P\otimes  TP)\xrightarrow{Id\otimes\ell_{k+1}}\wedge^{q-1}T^*P\otimes(\Delta_{k,~k+1}\hspace{1  mm}\rightarrow\hspace{1 mm} 0
\end{equation*}

\begin{equation*}
\hspace{15 mm}\downarrow\delta\hspace{43 mm}\downarrow\delta\hspace{54 mm}\downarrow\delta\hspace{10 mm}
\end{equation*}

\begin{equation*}
\hspace{2 mm}0\hspace{2 mm}\rightarrow\hspace{2 mm}\wedge^q T^*P\otimes\mathfrak{g}_{\ell+k}\hspace{2 mm}\longrightarrow\hspace{2 mm}\wedge^q T^*P\otimes(S^{\ell+k}T^*P\otimes TP)\hspace{6 mm}\xrightarrow{Id\otimes\ell_k}\hspace{5 mm}\wedge^q T^*P\otimes \Delta_{k-1,~k}\hspace{1 mm}\longrightarrow\hspace{1 mm} 0
\end{equation*}

\begin{equation}
\hspace{15 mm}\downarrow\delta\hspace{43 mm}\downarrow\delta\hspace{54 mm}\downarrow\delta\hspace{10 mm}
\end{equation}

\begin{equation*}
\hspace{2 mm}0\rightarrow\wedge^{q+1}T^*P\otimes\mathfrak{g}_{\ell+k-1}\rightarrow\wedge^{q+1}T^*P\otimes(S^{\ell+k-1}T^*P\otimes TP)\xrightarrow{Id\otimes\ell_{k-1}}\wedge^{q+1}T^*P\otimes \Delta_{k-2,~k-1}\rightarrow 0
\end{equation*}

\begin{equation*}
\hspace{15 mm}\downarrow\delta\hspace{43 mm}\downarrow\delta\hspace{70 mm}
\end{equation*}

\begin{equation*}
0\rightarrow\wedge^{q+2}T^*P\otimes\mathfrak{g}_{\ell+k-2}\rightarrow\wedge^{q+2}T^*P\otimes(S^{\ell+k-2}T^*P\otimes TP)\hspace{60 mm}
\end{equation*}
\end{sidewaysfigure}

\vspace{4 mm}

\newtheorem{faker}[LemmaCounter]{Lemma}
\begin{faker}
Let $Z\in J_{\mu}E$ be an element such that the first horizontal sequence in $(7)$ is surjective. Setting $h_{\mu}=(\Delta_{\mu-1,\mu})$ and indicating by $h_k,~k\geq\mu$ the $(k-\mu)-th$ algebraic prolongation of $h_{\mu}$, then  $h_k$ is $p-$acyclic if and only if $\mathfrak{g}_{\ell+k-1}$ is $(p+1)-$acyclic.
\end{faker}

\vspace{4 mm}

\newtheorem{fakers}[TheoremCounter]{Theorem}
\begin{fakers}
The hypotheses being those of the Lemma $14$, let us further assume that there exists a formal solution $(Z_k)_{k\geq\mu}$ of the fundamental system  $\textbf{F}=(\textbf{F}_k)_{k\geq\mu}$ (i.e., $Z_k\in\textbf{F}_k$ and $\rho_k Z_{k+1}=Z_k$) such that the first sequence of $(7)$ is surjective at the 
end along $(Z_k)$. This being so, we infer that

\vspace{4 mm}

1) the same property of surjectivity also holds for any other formal solution $W=(W_k)$ of $\textbf{F}$,

\vspace{4 mm}

2) $\textbf{F}_{k+1}=p\textbf{F}_k$ for $k\geq\mu_0$ and

\vspace{4 mm}

3) $\mu_0=aup~\{\mu,\mu_2\}$, where $\mu_2$ is the integer such that the symbol $\mathfrak{g}_{\ell+k},~k\geq\mu_2-1,$ becomes $2-acyclic$.

\vspace{4 mm}

\noindent
We recall that the symbol $\mathfrak{g}_{\ell+k}(Z_k),~Z_k\in\textbf{F}_k$, is entirely determined by $\mathfrak{g}_{\ell}(Z_0)$ and that the homological properties of this symbol are uniform along a trajectory $\textbf{F}_k$ in E (cf., the Proposition 4 in \cite{Kumpera2014}, sect. 6 and the Corollary 1). 
\end{fakers}

\vspace{4 mm}

\noindent
$\bf{Remark.}$ Under the hypotheses of the preceding Theorem, we have assumed the surjectivity of the 
morphisms $\ell_k$, for $k\geq\mu$, in view of being able to guarantee the stability of the symbols tangent 
to the fundamental equations $\bf{F_k}$ by means of the $2-acyclicity$ of the algebraic symbols $\mathfrak{g}_{\ell+k}~$. If we let down this property, very useful in applications, it will suffice to assume the surjectivity of the morphism $\ell_k$ for $k\geq\mu_0+1$.

\vspace{4 mm}

\newtheorem{fatal}[CorollaryCounter]{Corollary}
\begin{fatal}
The hypotheses being those of the last Theorem, the solutions of $\textbf{F}_{\mu_0}$ coincide with the simultaneous solutions of the (infinitely many) fundamental equations $\textbf{F}=(\textbf{F}_k)_{k\geq\mu}$ i.e., with the local sections S of E satisfying the conditions $im~j_k S\subset\textbf{F}_k,~k\geq\mu~$.
\end{fatal}

\vspace{4 mm}

\newtheorem{fractal}[CorollaryCounter]{Corollary}
\begin{fractal}
Let $\textbf{F}=(\textbf{F}_k)_{k\geq 0}$ be a family of integrable fundamental equations verifying the property $\rho_k\textbf{F}_{k+1}=\textbf{F}_k$ and let us assume further that $\ell_k,~k\geq\mu_0+1,$ is surjective along a formal solution 
of $\textbf{F},~\mu_0$ being the integer that stabilizes the symbols tangent to the equations $\textbf{F}_k~$. Under these conditions, we claim that $\textbf{F}_{k+1}=p\textbf{F}_k$ for all $k\geq\mu_0$ and that the solutions of $\textbf{F}_{\mu_0}$ coincide 
with the simultaneous solutions of the system  $\textbf{F}$. We also conclude thereafter that any particular solution of a given equation $\textbf{F}_k$ is as well a simultaneous solution. 
\end{fractal}

\vspace{4 mm}

The \textit{sorites} \footnote{\textbf{Sorite.} \textit{Du grecque "sorkites". Argument composé d'une suite de propositions liées entre elles de manière que l'attribut de chacune d'elles devienne le sujet de la suivante, et ainsi de suite, jusqu'à la conclusion, qui a pour sujet le sujet de la première et pour attribut 
l'attribut de l'avant dernière.} Larousse, \textit{quoque turbatio.}} being terminated (17 pages), let us 
return to geometry and try to do some adequate work. We thus begin by considering a finite prolongation space 
$(E,\pi,P,p)$ of order $\ell$ and take a local or global section \textit{S} of \textit{E} namely, a 
structure of species \textit{E}. We shall say that \textit{S} is \textit{homogeneous} or \textit{transitive 
in the base space} if the germs of \textit{S} at any two points $y,~y'\in\alpha(S)$ (the source of \textit{S}) are always equivalent. In other terms, there exists a germ of local diffeomorphism $\phi$ of \textit{P} with source \textit{y} and target \textit{y'} such that

\begin{equation*}
\phi(\underline{S}_y)\stackrel{def}{=}p\phi\circ\underline{S}_y\circ\phi^{-1}=\underline{S}_{y'}
\end{equation*}

\vspace{4 mm}

\noindent
where $\underline{S}_y$ denotes the germ of \textit{S} at the point \textit{y}. Inasmuch, we shall say that \textit{S} is \textit{homogeneous of order k} when the $k-$jets of \textit{S} at any two arbitrary points $y,~y'\in\alpha(S)$ are $k-th$ order equivalent which means that there exists an invertible jet $Y\in\Pi_{\ell+k}P$ such that $Y(j_kS(y))=j_kS(y'),~y=\alpha(Y),~y'=\beta(Y)$, and where, by definition, $Y(j_kS(y))=Y^{(k)}\cdot j_kS(y)\cdot Y^{-1}, Y^{(k)}$ being the invertible $k-$jet on \textit{E} that corresponds, by prolongation, to the $(\ell+k)-$jet \textit{Y} on \textit{P}. We shall say, finally, that \textit{S} is \textit{formally homogeneous} when the above $k-th$ order condition is verified for all integer \textit{k}.

\vspace{4 mm}

\noindent
We now assume that \textit{S} is homogeneous of order \textit{k} which means that the set $\{j_k S(y)|y\in\alpha(S)\}$ is contained in a single trajectory, denoted by $\Omega_k(S)$, under the action of the pseudo-group $\Gamma_{ell+k}$ on $J_k E$, natural prolongation to $k-th$ order of the \textit{general} pseudo-group $\Gamma(P)$ of all local diffeomorphisms of the manifold \textit{P}. Furthermore, when $\alpha(S)$ ia connected, the above set $\Omega_k(S)$ is also contained in a single trajectory $\omega_k(S)$ of the infinitesimal action $\mathcal{L}_{\ell+k}$ since these last trajectories are the connected components of the finite trajectories of $\Gamma_{ell+k}$. Each $\Omega_k(S)$ is a finite fundamental equation of order \textit{k} and each  $\omega_k(S)$ is an infinitesimal fundamental equation of the same order. When \textit{S} is formally homogeneous, we shall denote by $\Omega(S)=(\Omega_k(S))$ the family of all the respective fundamental equations. If, further, $j_k D$ takes all its values in a given infinitesimal trajectory, we denote by  $\omega(S)=(\omega_k(S))$ the corresponding family.

\vspace{4 mm}

\newtheorem{fatality}[PropositionCounter]{Proposition}
\begin{fatality}
    Let S be a homogeneous structure of order k and $\Omega_k(S)$ its finite fundamental equation of the same order. Then $\Omega_{h+1}(S)\subset p\Omega_h(S)$ for every integer $h<k$. If furthermore S is formally homogeneous, then there exists an integer $\mu_0$ such that $\Omega_{h+1}(S)$ is an open sub-bundle of $p\Omega_h(S)$ for any $h\geq\mu_0$. The integer $\mu_0$ is the order of stability ($1-$acyclicity) of the Spencer $\delta-$complex associated to the symbols tangent to the equations $\Omega_k(S)$. If, moreover, the symbol (i.e., the fibre) of $\Omega_{h+1}(S)$ above $\Omega_h(S)$ is, along $j_{\infty}S(y)$ and for $h\geq\mu$, an affine sub-space of the total symbol (more generally, if each connected component of the symbol is an affine sub-space), then $\Omega_{h+1}(S)=p\Omega_h(S)$ for $h\geq h_0=suo~\{\mu_0,\mu\}$. This being so, the solutions of $\Omega_{h_0}(S)$ coincide with the simultaneous solutions of $\Omega(S)$. The same properties also remain valid for $\omega_{h+1}(S)$ and $\omega(S)$, the integers $\mu_0$ and $h_0$ being the same as those for $\Omega(S)$.
\end{fatality}

\vspace{4 mm}

\noindent
It is clear that the solutions $S'$ of $\Omega_k(S)$ are precisely the structures of species \textit{E} $k-th$ order equivalent to the \textit{model S}, the $k-th$ order equivalence taking place for all couples of points $y\in\alpha(S)$ and $y'\in\alpha(S')$. The structures \textit{S'} are of course homogeneous of order \textit{k} and any structure of this type whose jet at a single point is $k-th$ order equivalent to a $k-$jet of \textit{S} is a structure $k-th$ order equivalent to \textit{S}. Likewise, the simultaneous solutions of the family $\Omega(S)$ are the structure formally equivalent to the model \textit{S}, such structures being all formally homogeneous. We shall say that a structure is \textit{connected} whenever $\alpha(S)$ is connected, a \textit{connected component} of an arbitrary structure \textit{S} being, by definition, the restriction of \textit{S} to a connected component of $\alpha(S)$. Its image is a connected component of \textit{im S}.

\vspace{4 mm}

\newtheorem{reality}[PropositionCounter]{Proposition}
\begin{reality}
Let S be a formally homogeneous structure. Then, the local or global structures of species E and formally equivalent to the model S are the solutions of $\Omega(S)$.  Moreover, when the hypotheses of the last proposition concerning the affine nature of the symbols of $\Omega-h(S)$ are verified, the formal equivalence is then a consequence of the equivalence at order $h_0$. Finally, when $\omega(S)$ is defined (which occurs especially when S is connected), the connected structures formally equivalent to the model S are the solutions of $\omega(S)$.
\end{reality}

\vspace{4 mm}

\noindent
The system $\Omega(S)$ is, according to the terminology of \cite{Molino1972}, \textit{the fundamental differential system} for the structures of species \textit{E} that are formally equivalent to \textit{S}. In his article, the author considers the model structure defined on a manifold $P_0$ eventually distinct from \textit{P}. In ou case, we simply identify $P_0$ to an open set of \textit{P} since it is always possible to \textit{transfer} the model given on $P_0$ onto a model defined on an open set of \textit{P} without, for that matter, modifying the equivalence relation.

\vspace{4 mm}

\noindent
We could hereafter contemplate in defining the \textit{reduced} fundamental system and in transcribing some of the results of \cite{Molino1972}. However, we leave such transcriptions, not all together evident, to the care of the reader since, at present, we are inclined to examine other important aspects of the theory.

\vspace{4 mm}

\noindent
Let us first observe that the different notions of homogeneity introduced above refer essentially to the pseudo-group $\Gamma(S)$ of all automorphisms of \textit{S} as well as to the groupoids $\mathcal{R}_{\ell+k}(S)$. In fact, \textit{S} is \textit{homogeneous} when $\Gamma(S)$ is transitive, \textit{homogeneous of order k} when the groupoid $\mathcal{R}_{\ell+k}(S)$ is transitive and finally \textit{formally homogeneous} when the transitivity of $\mathcal{R}_{\ell+k}(S)$ occurs for all \textit{k}. We can therefore introduce the corresponding infinitesimal notions namely, we can say that \textit{S} is \textit{infinitesimally homogeneous} when $\mathcal{L}(S)$ is transitive, $k-th$ \textit{order infinitesimally homogeneous} when $\textbf{R}_{\ell+k}(S)$ is transitive i.e., $\beta(\textbf{R}_{\ell+k}(S))_y=T_yP$, for all $y\in P$, and \textit{formally infinitesimally homogeneous} when the last condition is verified for all \textit{k}. When \textit{S} is only defined on an open set $\mathcal{U}$, we shall simply replace, in the above definitions, the manifold \textit{P} by the open set $\mathcal{U}$.

\vspace{4 mm}

\noindent
Let \textit{S} be a structure of species \textit{E} and let us inquire for conditions rendering the target map $\beta:(\textbf{R}_{\ell+k}S)_y~\longrightarrow~T_y P$ surjective. For this, we set $Y=j_{\ell+k}S(y)$ and recall that

\begin{equation*}
(\textbf{R}_{\ell+k}S)_y=\{j_{\ell+k}\xi(y)\in J_{\ell+k}TP~|~(p_k\xi)_Y\in T_Y(im~j_k S)\}.
\end{equation*}

\vspace{4 mm}

\noindent
However, since $\xi_y=\beta(j_{\ell+k}\xi(y))=\alpha_*(p_k\xi)_Y,~\alpha:J_kE\longrightarrow P$, we see promptly that the surjectivity condition is given by $T_Y(im~j_k S)\subset(\Delta_k)_Y$. Consequently, a structure \textit{S} is infinitesimally homogeneous of order \textit{k} if and only if $im~j_k S$ is an integral sub-manifold of the distribution (differential system) $\Delta_k$. Taking into account the integrability (involutivity) of $\Delta_k$ (\cite{Kumpera2014}, Theorem 2) and the fact that the integral leaves of $\Delta_k$ are the connected components of the trajectories $\Omega_k$ of $\Gamma_{\ell+k}$, we obtain the following result:

\vspace{4 mm}

\newtheorem{non-reality}[PropositionCounter]{Proposition}
\begin{non-reality}
The following statements are equivalent:

\vspace{4 mm}

1. S is infinitesimally homogeneous of order k.

\vspace{4 mm}

2. Each connected component of $~im~j_k S~$ is contained in an infinitesimal trajectory of $\mathcal{L}_{\ell+k}$.

\vspace{4 mm}

3. Each connected component of $~im~j_k S~$ is contained in a finite trajectory of $\Gamma_{\ell+k}$.
 
 \vspace{4 mm}
 
 4. Each connected component of $im~S$ is homogeneous of order k.
\end{non-reality}

\vspace{4 mm}

\newtheorem{non-real}[CorollaryCounter]{Corollary}
\begin{non-real}
The following statements hold:

\vspace{4 mm}

1. Every connected structure that is infinitesimally homogeneous of order k is also homogeneous of order k and every homogeneous structure of order k is also infinitesimally homogeneous of the same order.

\vspace{4 mm}

2. Let S be a connected and infinitesimally homogeneous model structure of order k. Then the class of all the connected structures S' equivalent to order k to the model S is given by the connected solutions of the fundamental equation $\omega_k(S)$. Each such equivalent structure S' is infinitesimally homogeneous of order k.

\vspace{4 mm}

3.  Let S be a connected and infinitesimally formally homogeneous model structure.  Then the class of all the connected structures S' formally equivalent to the model S is given by the connected solutions of the infinitesimal fundamental system $\omega(S)$ or, when the hypotheses of the Proposition 1 are verified, by the connected solutions of $\omega_{h_0}(S)$ in which case the formal equivalence is a consequence of the equivalence at order $h_0$. Any solution S' is infinitesimally formally homogeneous.
\end{non-real}

\vspace{4 mm}

\noindent
We next discuss the notion of \textit{transitivity}. The structure \textit{S} is said to be \textit{transitive of order k} or $k-th$ order \textit{transitive} when it is homogeneous of order \textit{k} and if further the projections $\mathcal{R}_{\ell+k}(S)~\longrightarrow~\mathcal{R}_{\ell+h}(S),~0\leq h\leq k~,$ are all surjective (in fact, surmersions on account of the Lemma 8). The structure \textit{S} is said \textit{formally transitive} when the above conditions are verified for all \textit{k} and, finally, just \textit{transitive} when it is formally transitive and, moreover, when $\mathcal{R}_{\ell+k}(S)=J_{\ell+k}\Gamma(S)$. In other terms, this means that any equivalence of finite order among the jets of \textit{S} is actually realised (achieved) by the jet of a local automorphism of \textit{S}. In essentially the same way, we can introduce the notions of $k-th$ order infinitesimal transitivity, formal infinitesimal transitivity and, for short, transitivity by simply replacing $\Gamma(S)$ by $\mathcal{L}(S)$ and $\mathcal{R}_{\ell+k}(S)$ by $\textbf{R}_{\ell+k}(S)$.

\vspace{2 mm}

\noindent
We observe that the notions of $k-th$ order and formal finite or infinitesimal transitivity, are invariant by $k-th$ order and formal equivalences respectively. Thus, if we start with a formally transitive resp., $k-th$ order model, every solution of the fundamental system $\Omega(S)$ resp, $\Omega_k(S)$ is a formally transitive resp., $k-th$ order structure. The same remarks apply of course for the infinitesimal context taking, however, the care in taking for \textit{S} a formally homogeneous and infinitesimally formally transitive resp., $k-th$ order infinitesimally transitive structure.

\vspace{2 mm}

\noindent
Let us now make abstracrion of the model. It is clear that a structure \textit{S} is homogeneous of order \textit{k} if and only if this structure is a solution of a fundamental equation $\Omega_k$. Moreover, each connected component of \textit{S} is a solution of an infinitesimal equation $\omega_k$ contained in $\Omega_k$. Inasmuch, \textit{S} is infinitesimally homogeneous of order \textit{k} if and only if each connected component is a solution of a finite or infinitesimal fundamental equation of order \textit{k}, though these equations may vary along with the connected component. We have analogous conclusions in the formal case, the number of equations being infinite unless the hypotheses of the Proposition 1 be verified.

\vspace{2 mm}

\noindent
We next inquire on the equations bearing on the finite jets of a structure \textit{S} and displaying their $k-th$ order or their formal transitivity and start examining the infinitesimal aspect that will curiously place in evidence a new element namely, the \textit{Medolaghi-Vessiot equations}.

\vspace{2 mm}

\noindent
Let us be given a structure \textit{S} of species \textit{E}, let us fix an integer \textit{k} and let us pose ourselves the following problem. What are the conditions to which are bound the finite jets $j_{k+1}S(y)~,~y\in\alpha(S)\subset P$ in such a way that

\vspace{2 mm}

\textit{a}) the projection $\textbf{R}_{\ell+k+1}(S)~\longrightarrow~\textbf{R}_{\ell+k}(S)$ be surjective ($1-$integrable at order $\ell+k$) and

\vspace{2 mm}

\textit{b}) $\textbf{R}_{\ell+k}(S)~\longrightarrow~TP$ be surjective (infinitesimal homogeneity of order \textit{k}).

\vspace{2 mm}

\noindent
We observe immediately that the second condition can be replace, in view of (a), by the surjectivity of $\textbf{R}_{\ell+k+1}(S)~\longrightarrow~TP$. Moreover, we can replace the above problem by the equivalent problem

\vspace{2 mm}

$a_1)~\textbf{R}^0_{\ell+k+1}(S)~\longrightarrow~\textbf{R}^0_{\ell+k}(S)$ be surjective and

\vspace{2 mm}

$b_1)~\textbf{R}_{\ell+k+1}(S)~\longrightarrow~TP$ be surjective,

\vspace{2 mm}

\noindent
since the \textit{small} diagram hereafter is commutative and exact.

\begin{equation*}
0\longrightarrow\textbf{R}^0_{\ell+k+1}(S)\longrightarrow\textbf{R}_{\ell+k+1}(S)\longrightarrow TP\longrightarrow 0
\end{equation*}
\begin{equation*}
\hspace{7 mm}\downarrow\hspace{25 mm}\downarrow\hspace{18 mm} |\!|
\end{equation*}
\begin{equation*}
0\longrightarrow\hspace{3 mm}\textbf{R}^0_{\ell+k}(S)\hspace{2 mm}\longrightarrow\hspace{2 mm}\textbf{R}_{\ell+k}(S)\hspace{1 mm}\longrightarrow TP\longrightarrow 0
\end{equation*}
\begin{equation*}
\hspace{12 mm}\vdots\hspace{26 mm}\vdots\hspace{27 mm}
\end{equation*}
\begin{equation*}
\hspace{13 mm} 0\hspace{25 mm} 0\hspace{27 mm}
\end{equation*}

\vspace{2 mm}

\noindent
We first examine the problem $a_1$. The relations (3) lead us to introduce the equation $\Theta_{k+1}$ of order $k+1$ defined on the fibration $E~\longrightarrow~P$

\begin{equation*}
\Theta_{k+1}=\{Z\in J_{k+1}E~|~(\textbf{R}^0_{\ell+k+1})_Z\longrightarrow (\textbf{R}^0_{\ell+k})_{\rho(Z}\longrightarrow 0\hspace{2 mm}is~exact\}.
\end{equation*}

\vspace{2 mm}

\noindent
The following Lemma is evident on account of the Lemma 1.

\vspace{4 mm}

\newtheorem{non-complex}[LemmaCounter]{Lemma}
\begin{non-complex}
The equation $\Theta_{k+1}$ is invariant under the finite action of $\Gamma_{\ell+k+1}$ as well as under the infinitesimal action of $\mathcal{L}_{\ell+k+1}$ or, in other terms, $\Theta_{k+1}$ is a union of trajectories. Inasmuch, this equation is invariant under the left action of the groupoid $\Pi_{\ell+k+1}P$ on $J_kE$ and the infinitesimal action of the algebroid $J_{\ell+k+1}TE$. The solutions of $\Theta_{k+1}$ are the local or global structures of species E verifying the condition $a_1$.     
\end{non-complex}

\vspace{2 mm}

\noindent
Since the condition $b_1$ is no other than the $(k+1)-$order infinitesimal homogeneity condition, we derive the following result:

\vspace{4 mm}

\newtheorem{complex}[PropositionCounter]{Proposition}
\begin{complex}
A structure of species E verifies the conditions $(a)$ and $(b)$ if and only if each of its connected components is a solution of a finite or infinitesimal fundamental equation contained in the equation $\Theta_{k+1}$. The connected structures satisfying $(a)$ and $(b)$ are the connected solutions of the infinitesimal dundamental equations contained in $\Theta_{k+1}$.
\end{complex}

\vspace{2 mm}

\noindent
We thus see that there exists a whole family of differential equations of order \textit{k+1} on the fibration $E~\longrightarrow~P$ solving the proposed problem. Whereas the fundamental equations are well behaved with respect to prolongations, this does not happen with the equations $\Theta_{k+1}$. In fact, there does not exist, \textit{a priori}, any plausible relation between $\Theta_{k+1}$ and $p\Theta_k$. Inasmuch, we can even note that there does not exist any apparent relation between $\rho_k(\Theta_{k+1})$ and  $\Theta_k$, all depending on the specific structure of the prolongation space \textit{E}. We shall therefore search for the desired results with the help of the Theorem of Quillen.

\vspace{2 mm}

\noindent
With this in mind, let us first recall that this theorem strongly relies on the $2-$acyclicity of the symbols and the Corollary 1 tells us that the homological properties of the symbols $\mathfrak{g}_{\ell+k}$ are uniform along the inverse image of an orbit of degree zero hence, in particular, along a fundamental equation $\bf{F}_k$. On the other hand, if \textit{S} is a solution of $\bf{F}_k$, \textit{S} takes its values in the orbit of order zero $\textbf{F}_0=\beta(\textbf{F}_k)$. Consequently, when $\bf{F}_k$ is integrable then $\textbf{F}_k\subset J_k\textbf{F}_0$ and this leads us naturally to only consider the prolongation space $(\textbf{F}_0,\pi,P,p)$ of order $\ell$ that will be called \textit{subordinate} to $(E,\pi,P,p)$ and whose prolongation structure is simply obtained by restricting the structure of \textit{E} (\textit{cf.} \cite{Kumpera2014}, Theorem 2). When $\bf{F}_0$ is a finite trajectory, the prolongation laws of \textit{E} admit natural restrictions. Quite to the contrary, when the trajectory is infinitesimal, the infinitesimal prolongation law restricts without any problem while the finite law will have to be restricted, at each order \textit{k}, to the open sub-groupoid of $\Pi_{\ell+k}P$, composed of those finite jets that preserve the $k-th$ order transverse contact elements of dimension \textit{n} tangent to $\bf{F}_0$ but also restricted to the solutions of this sub-groupoid. In particular, the $\alpha-$connected component $(\Pi_{\ell+k}P)_0$ is contained in this sub-groupoid. It is clear that $J_k\textbf{F}_0$ is invariant in $J_kE$ (always respecting the above restrictions concerning the infinitesimal orbits) and we shall denote by $\Theta_{k+1}(\textbf{F}_0)$ and $\mathfrak{g}_{\ell+k}(\textbf{F}_0$ the equation $\Theta_{k+1}$ and the symbol $\mathfrak{g}_{\ell+k}$ respectively, when restricted to $J_{k+1}\textbf{F}_0$ and $J_k\textbf{F}_0~$. The restricted symbol has the same homological properties everywhere. Lastly, we shall call $\textbf{F}_0-admissible$ any fundamental equation contained in $\Theta_{k+1}$ thus finding the so claimed Medolaghi-Vessiot equations (\cite{Vessiot1903}, p.436, eqs.(58) and (59)).

\vspace{4 mm}

\newtheorem{sub-complex}[TheoremCounter]{Theorem}
\begin{sub-complex}
Let $(E,\pi,P,p)$ be a finite prolongation space of order $\ell~$, $(\textbf{F}_0,\pi,P,p)$ the prolongation space subordinate to an orbit $\textbf{F}_0$ and $\eta_2=\eta_2(\textbf{F}_0)$ the integer from whereon the symbol $\mathfrak{g}_{\ell+k}(\textbf{F}_0)~,~k\geq\eta_2~$, becomes $2-$acyclic.The linear equation $\textbf{R}_{\ell+k}(S)$ associated to any solution of an $\textbf{F}_0-$admissible fundamental equation of order $k+1\geq\eta_2+1$ is then transitive and formally integrable \textit{i.e.}, formally transitive. Moreover, when $\eta_{\infty}=\eta_{\infty}(\textbf{F}_0)$ is the integer where from the restricted symbol $\mathfrak{g}_{\ell+k}(\textbf{F}_0)~,~k\geq\eta_{\infty}~$, becomes involutive, then the equation $\textbf{R}_{\ell+k}(S)~,~k\geq\eta_{\infty}~$, is involutive.
\end{sub-complex}

\vspace{4 mm}

\noindent
The argument is as follows. \textit{S} being a solution of a fundamental equation of order $k+1$, we know (\textit{cf.} the Lemma 8) that $\textbf{R}_{\ell+k}(S)$ as well as $\textbf{R}_{\ell+k+1}(S)$ are regular equations and that  $\textbf{R}_{\ell+k+1}(S)=p\textbf{R}_{\ell+k}(S)$. Furthermore, the symbol $\mathfrak{g}_{\ell+k}(S)=(j_kS)^{-1} \mathfrak{g}_{\ell+k}(\textbf{F}_0)$ is $2-$acyclic and the morphism $\textbf{R}_{\ell+k+1}(S)~\longrightarrow~ \textbf{R}_{\ell+k}(S)$ is surjective since \textit{S} is a solution of the Medolaghi-Vessiot equations of order $k+1$. We thus re-encounter the hypotheses of Quillen's Theorem (\cite{Kumpera1972}) and consequently $\textbf{R}_{\ell+k}(S)$ is formally integrable. We also observe, \textit{en passant}, that $\textbf{R}_{\ell+k+2}(S)=p\textbf{R}_{\ell+k+1}(S)$ since the equation of order $\ell+k+1$ is regular and the $2-$acyclicity implies that that $\textbf{R}_{\ell+k+2}(S)$ is also regular. An inductive argument will then show that all the equations $\textbf{R}_{\ell+k+h}(S)$ are regular and that $\textbf{R}_{\ell+k+h+1}(S)=p\textbf{R}_{\ell+k+h}(S)$, which by the way is already stated in the Theorem of Quillen. We finally observe that this Theorem does not require the regularity of the first two equations since it suffices in fact that $\textbf{R}_{\ell+k}(S)$, hence also its prolongation, be defined by a differential operator which is always the case for the equations $\textbf{R}_{\ell+k}(S)$. The regularity is just a Black Friday extra bonus.

\vspace{4 mm}

\newtheorem{no-faker}[CorollaryCounter]{Corollary}
\begin{no-faker}
When $\mathfrak{g}_{\ell}(\textbf{F}_0)$ is $2-$acyclic then every solution S of an $\textbf{F}_0-$admissible fundamental equation of order $1$ is an infinitesimally formally transitive structure of species E (more so, of species $\textbf{F}_0$). If, further, $\mathfrak{g}_{\ell}(\textbf{F}_0)$ is involutive, the same will be true for the structure S i.e., for the equation $\textbf{R}_{\ell}(S)$.  
\end{no-faker}

\vspace{4 mm}

\noindent
We can now consider the section $j_kS$ as being a structure of species $J_kE$, \textit{prolongation of order k} of the structure \textit{S} and the last Theorem can be paraphrased by the

\vspace{4 mm}

\newtheorem{substantial}[CorollaryCounter]{Corollary}
\begin{substantial}
Let $(\textbf{F}_0,\pi,P,p)$ be the prolongation space subordined to an orbit $\textbf{F}_0$ of $(E,\pi,P,p)$. The $k-th$ order prolongation of any $\textbf{F}_0-$admissible fundamental equation of order $k+1\geq\eta_2+1~$ is a structure of species $J_kE$ (viz. of species $J_k$) that moreover is infinitesimally formally transitive. Inasmuch, the prolongation of order $\geq\eta_{\infty}$ is an involutive structure 
\end{substantial}

\vspace{4 mm}

\noindent
We need not play, for the prolonged structures, the Ehresmannian semi- or non-holonomic game \textit{i.e.}, by considering spaces like $J_h(J_{\ell+k}E)$, since we can see promptly that $\textbf{R}_{(\ell+k)+h}(j_k S)=\textbf{R}_{\ell+k+h}(S)$ as soon as we define these equations by the \textit{contact order} of the prolonged vector fields with the sections defining respectively the structure.

\vspace{4 mm}

\newtheorem{substantially}[CorollaryCounter]{Corollary}
\begin{substantially}
Let $\textbf{F}_{k+1}$ be an integrable fundamental equation contained in $\Theta_{k+1}$ and assume that $\mathfrak{g}_{\ell+k}(\textbf{F}_0)~,~\textbf{F}_0=\rho_0\textbf{F}_{k+1}~$, is $2-$acyclic. Then, for any $h\geq 0~$, the $h-th$ prolongation $p_h\textbf{F}_{k+1}$ is an $\textbf{F}_0-$admissible fundamental equation of order $k+h+1$ namely, the unique integrable fundamental equation projecting onto $\textbf{F}_{k+1}$.
\end{substantially}

\vspace{4 mm}

\noindent

\vspace{4 mm}

a) The structure \textit{S} is formally homogeneous hence $~im~j_{k+h}S\subset\textbf{F}_{k+h}~$, $i.e.,~\textbf{F}_{k+h}$ is integrable and consequently

\begin{equation*}
\textbf{F}_{k+h}\subset p\textbf{F}_{k+h}\subset p_h\textbf{F}_{k+1}~.
\end{equation*}

\vspace{2 mm}

b) The morphism

\begin{equation*}
(\textbf{R}^0_{\ell+k+h})_{Y_{k+h}}\longrightarrow (\textbf{R}^0_{\ell+k+h-1})_{Y_{k+h-1}}
\end{equation*}

\vspace{2 mm}

\noindent
is surjective hence $\textbf{F}_{k+h}$ is $\textbf{F}_0-$admissible and, on account of the Lemma 12, the mapping $(\ell_{k+h})_{Y_{k+h}}$ is also surjective. The $2-$acyclicity of $\mathfrak{g}_{\ell+k}(\textbf{F}_0)$ implies (Lemma 15) that $\mu_0=k+1$ is equal to the integer that stabilizes the symbols tangent to the equations $\textbf{F}_k$ and consequently (Corollary 10), $\textbf{F}_{k+h+1}=p\textbf{F}_{k+h}$ or, in other terms, $\textbf{F}_{k+h+1}=p_h\textbf{F}_{k+1}~$. The uniqueness of the equation $\textbf{F}_{k+h+1}$ is obvious since any other integrable fundamental equation of order $k+h+1$ that projects onto $\textbf{F}_{k+1}$ is necessarily contained hence equal to $p_h\textbf{F}_{k+1}~$.

\vspace{4 mm}

\noindent
This Corollary shows that given an $\textbf{F}_0-$admissible and integrable fundamental equation $\textbf{F}_{k+1}~,~k+1\geq\eta_2+1$ (\textit{cf.} the Theorem 3), there will be no place in pushing any further the calculations that is to say, search for solutions \textit{S} of $\textbf{F}_0-$admissible fundamental equations $\textbf{F}_{k+1+h}$ projecting upon $\textbf{F}_{k+1}$ since we shall find no other than the solutions of $\textbf{F}_{k+1}$, no additional restriction being therefore possible.

\vspace{4 mm}

\noindent
We thus see that the local or global structure of species \textit{E} admitting an \textit{infinitesimally formally transitive} prolongation can be searched for among the solutions of rhe $\textbf{F}_0-$admissible and integrable fundamental equations of orders $k+1\geq\eta_2(\textbf{F}_0)+1$. Moreover, this method exhausts the connected structures. In fact, if $j_k S$ is an infinitesimally formally transitive structure, the equation $\textbf{R}_{\ell+k}(S)$ is formally integrable and transitive. Taking, if necessary, an integer $k'>k~$, this equation acquires a $2-$acyclic symbol. Since, by definition, this equation verifies the conditions $(a)$ and $(b)$, the section \textit{S} is the solution of a fundamental equation $\textbf{F}_{k+1}$ contained in $\Theta_{k+1}$ (Proposition 3) and thereafter is  $\textbf{F}_0-$admissible, where  $\textbf{F}_0=\rho_0\textbf{F}_{k+1}~$. Since $\mathfrak{g}_{\ell+k}(\textbf{F}_0)$ is $2-$acyclic, we can find the structure \textit{S} by the method of the Theorem 3. When the structure is not connected, each of its components can be determined by the above method.

\vspace{4 mm}

\noindent
The preceding results involve integrable fundamental equations hence it is high time for us to examine the formal integrability of the less restrictive  $\textbf{F}_0-$admissible fundamental equations. The next Lemma as well as the Theorem show that, subject to reasonable hypotheses, the Medolaghi-Vessiot equations have good properties with respect to prolongations and formal integrability.

\vspace{4 mm}

\newtheorem{insubstantial}[LemmaCounter]{Lemma}
\begin{insubstantial}
Let $\textbf{F}_{k+1}$ be an $\textbf{F}_0-$admissible fundamental equation verifying the two properties: $p\textbf{F}_{k+1}\neq\phi$ and $\mathfrak{g}_{\ell+k}(\textbf{F}_0)$ is $2-$acyclic. Then $p\textbf{F}_{k+1}$ is also an $\textbf{F}_0-$admissible fundamental equation.  
\end{insubstantial}

\vspace{4 mm}

\noindent
In order to prove this rather remarkable result, we first observe that the condition $p\textbf{F}_{k+1}\neq\phi$ entails (Corollary 2) that $p\textbf{F}_{k+1}$ is a regular equation of order $k+2$ namely, a locally trivial affine sub-bundle of $J_{k+2}E\longrightarrow\textbf{F}_{k+1}$ and, in particular, that its tangent symbol $(\gamma_{k+2})_Z~$, at a point $Z\in p\textbf{F}_{k+1}$ is the algebraic prolongation of $(\Delta_{k,k+1})_Y~,~Y=\rho_{k+1}Z~$, the latter being the symbol tangent to $\textbf{F}_{k+1}$ at the point $Y~$. We also know (Lemma 11) that $p\textbf{F}_{k+1}$ is a union of trajectories \textit{i.e.,} the trajectory $\textbf{F}_{k+2}$ that contains the point $Z$ is entirely contained in $p\textbf{F}_{k+1}~$. We infer that the tangent symbol $(\Delta_{k+1,k+2})_Z$ of $\textbf{F}_{k+2}$ at the point $Z$ is contained in $(\gamma_{k+2})_Z~$. Let us now re-examine the diagram $(8)$ at its stage $(1,k+1)$ where we shall replace the kernel  $\Delta_{k-1,k}$ by $im~\ell_k~$. The fact that $\textbf{F}_{k+1}$ be $\textbf{F}_0-$admissible implies (Lemma 12) that $im~(\ell_{k+1})_Y=(\Delta_{k,k+1})_Y$ and the proper definition of $\ell_k$ (\textit{cf.}, $(7)_k$) shows that $im~(\ell_{k+2})_Z\subset (\Delta_{k+1,k+2})_Z~$. Last but not least, the $2-$acyclicity of $\mathfrak{g}_{\ell+k}(\textbf{F}_0)$ implies that the third vertical sequence in $(8)$ \textit{i.e.}, the sequence

\begin{equation*}
im~(\ell_{k+2})_Z\xrightarrow{\delta}T^*_y P~\otimes~(\Delta_{k,~k+1})_Y\xrightarrow{\delta}\wedge^2 T^*_y~\otimes~im~(\ell_k)_Y
\end{equation*}

\vspace{2 mm}

\noindent
$Y=\rho_k Z~$, is exact. But this means precisely that $im~(\ell_{k+2})_Z$ is the algebraic prolongation of $(\Delta_{k,k+1})_Y$ or, in other terms, that $im~(\ell_{k+2})_Z=(\gamma_{k+2})_Z\supset\Delta_{k+1,k+2})_Z$ and consequently that $im~(\ell_{k+2})_Z=\Delta_{k+1,k+2})_Z~$. The trajectory $\textbf{F}_{k+2}$ is $\textbf{F}_0-$admissible (Lemma 12) and  $\textbf{F}_{k+2}=p\textbf{F}_{k+1}$ since $\ell_{k+2})_Z$ is surjective (Corollary 6) and the symbols tangent to both of these equations coincide (\textit{cf.}, the proof of the Lemma 14).

\vspace{4 mm}

\newtheorem{disturbing}[TheoremCounter]{Theorem}
\begin{disturbing}
Let $\textbf{F}_{k+1}$ be an $\textbf{F}_0-$admissible fundamental equation and let us further assume that the symbol $\mathfrak{g}_{\ell+k}(\textbf{F}_0)$ is $3-$acyclic. Then $\textbf{F}_{k+1}$ is formally integrable if and only if $p\textbf{F}_{k+1}\neq\phi$ (on each fibre!) and, this being the case, each $p_h\textbf{F}_{k+1}$ is an $\textbf{F}_0-$admissible fundamental equation. If, moreover, 
$\mathfrak{g}_{\ell+k}(\textbf{F}_0)$ is involutive then $\textbf{F}_{k+1}$ is also involutive. 
\end{disturbing}

\vspace{4 mm}

\noindent
The proof is almost obvious. We know, according to the Corollary 2, that $p\textbf{F}_{k+1}\neq\phi$ if and only if $\textbf{F}_{k+1}$ is $1-$integrable and, whenever this condition is fulfilled, that $p\textbf{F}_{k+1}$ is also a regular equation. More precisely, $p\textbf{F}_{k+1}$ is a locally trivial affine sub-bunbdle of $J_{k+2}E~\longrightarrow~\textbf{F}_{k+1}~$. The $3-$acyclicity hypothesis on $\mathfrak{g}_{\ell+k}(\textbf{F}_0)$ entails, having in mind the Lemma 16, that the tangent symbol $\Delta_{k,k+1}$ of $\textbf{F}_{k+1}$ is everywhere $2-$acyclic. We can then complete the proof using a non-linear version of the Theorem of Quillen and showing the formal integrability of $\textbf{F}_{k+1}$ by induction on the integer \textit{h} so as to propagate the result of the preceding Lemma to all the prolongations (\textit{see, for a hint} \cite{Kumpera1972} and \cite{Kumpera1999}).

\vspace{4 mm}

\newtheorem{best}[]{Scholium}
\begin{best}
When $(E,\pi,P,p)$ is an analytic prolongation space and when $\textbf{F}_{k+1}$ is an  $\textbf{F}_0-$admissible $1-$integrable fundamental equation whose symbol $\mathfrak{g}_{\ell+k}(\textbf{F}_0)$ is $3-$acyclic, then $\textbf{F}_{k+1}$ is integrable, everyone of its element being the $(k+1)-$jet of an analytic solution S. Furthermore, the solutions of any such equation are analytic structures whose prolongations of order $k+1$ are infinitesimally transitive.
\end{best}

\vspace{4 mm}

\noindent
Prolonging if necessary the equation $\textbf{F}_{k+1}$, we find an equation with involutive symbol hence the result as a consequence of \textit{the Cartan-Kähler or the Cauchy-Kovalevskaya (or even perhaps the Monge-Ampère?)} Theorems. The last assertion follows from the above mentioned theorems and an eventual prolongation procedure of the formally integrable equation $\textbf{R}_{\ell+k}(S)$ so as to attain involutivity.

\vspace{4 mm}

\noindent
The condition $p\textbf{F}_{k+1}\neq\phi$ is a point wise imposition since, for this to be so, it is necessary and sufficient that, at a point $Y\in\textbf{F}_{k+1}~$, the tangent space $T_Y\textbf{F}_{k+1}$ contain a contact element of order 1 transverse (to the fibration $J_{k+1}E\longrightarrow P$) and \textit{holonomic}. This condition correspoonds, when $k=0~$, to what the author of \cite{Pomme1973} calls \textit{The generalized Jacobi relations}, that seem rather out of place.

\vspace{4 mm}

\noindent
Let us next show that the restrictions $\Theta_{k+1}(\textbf{F}_0)$ of $\Theta_{k+1}$ to the subordinate prolongation spaces are closed sub-sets of $J_{k+1}\textbf{F}_0~$. We shall prove simultaneously that the restrictions $\Theta_{k+1}\cap\beta^{-1}_{k+1}(\textbf{F}_0)$ are closed sub-spaces in $\beta^{-1}_{k+1}(\textbf{F}_0)$ where $\beta_{k+1}:J_{k+1}E\longrightarrow J_kE~$. The argument will be based on a local coordinates calculation which, in its essence, can already be found in \cite{Amaldi1942}, Vol.2, Chap.V, \S9 (\textit{Cenni sulle ricerche di Engel-Medolaghi-Vessiot}).

\vspace{4 mm}

\noindent
Let $(x^i)$ be a system of local coordinates defined on an open set \textit{U} of \textit{P}, $(x^i,y^{\lambda})$ an adapted system on the open set $\mathcal{U}$ of \textit{E} and $(x^i,y^{\lambda},y^{\lambda}_{\alpha})_{|\alpha|\leq k}$ the corresponding system on $\beta^{-1}_k(\mathcal{U})\subset J_k E~$. Let us further denote by $(z^i)$ the same coordinates $(x^i)$ on \textit{U} and then start our \textit{coordinate juggling game} that will be written in French so as to spare those readers that suffer intolerance with coordinates.

\vspace{4 mm}

\noindent
En chaque point fixé $z=(z^i)$ de \textit{U}, le champ de vecteurs  $\xi^{\beta}_{i,z}=(1/\beta!)(x-z)^{\beta}~\partial/\partial x^i~,~(x-z)^{\beta}=(x^1-z^1)^{\beta_1} \dots (x^n-z^n)^{\beta_n}~,~1\leq |\beta|\leq k~$, 
s'annule au point \textit{z} et par conséquent $\epsilon^{\beta}_i (z)=j_k~\xi^{\beta}_{i,z} (z)\in(J^0_k TP)_z~$. On remarque en plus que l'ensemble $\{\epsilon^{\beta}_i (z)\}$ est une base de cette algèbre de Lie, \textit{la base de Engel} suivant la terminologie de \cite{Pomme1973} car, lorsque $\xi=\sum~\xi^i~\partial /\partial x^i~$, alors $j_k\xi(z)=\sum~(\partial^{|\beta|}\xi^i/\partial x^{\beta})(z)\epsilon^{\beta}_i (z)~$.

\vspace{4 mm}

\noindent
Les sections $\epsilon^{\beta}_i :U~\longrightarrow~J^0_kTP$ sont différentiables et constituent una base locale du fibré $J^0_kTP|U~$. Plaçons nous maintenant dans le cadre d'un espace de prolongement \textit{E} d'ordre $\ell$ et reprenons les  notations ci-dessus avec \textit{k} remplacé par $\ell+k~$. À chaque section $\epsilon^{\beta}_i~,~1\leq |\beta|\leq \ell+k~$, de $J^0_kTP$ correspond un champ de vecteurs prolongé $p_k\epsilon^{\beta}_i$ de $J_kE$ défini sur l'ouvert $\alpha^{-1}_k(U)$ et prenant ses valeurs dans le fibré vectoriel $VJ_kE~$. L'ensemble de ces champs prolongés n'est pas en général linéairement indépendant en tous les points mais pourtant il engendre la distribution verticale $\Delta_k\cap VJ_kE~$. Or, pour tout $Z\in\alpha^{-1}_k(U)~$, l'isotropie $(\textbf{R}^0_{\ell+k})_Z$ est l'ensemble des combinaisons linéaires de la base de Engel d'ordre $\ell+k$ au point \textit{z} dont lea coefficients sont ceux qui établissent une relation de dépendance linéaire entre les vecteurs prolongéa $p_Z\epsilon^{\beta}_i~$. Pour mieux décrire l'espace vectoriel de ces relations linéaires revenons, avec les modifications évidentes, à la suite exacte (7) de \cite{Kumpera2014}

\begin{equation*}
0\longrightarrow \textbf{R}^0_{\ell+k}\longrightarrow J^0_{\ell+k}TP~\times_P~J_kE\xrightarrow{\overline{\lambda}_k} VJ_kE 
\end{equation*}

\vspace{2 mm}

\noindent
et posons

\begin{equation*}
\overline{\lambda}_k=\sum_{|\alpha|\leq k}~a^{(\lambda,\alpha)}~\partial /\partial y^{\lambda}_{\alpha}~.
\end{equation*}

\vspace{2 mm}

\noindent
Les coefficients $a^{(\lambda,\alpha)}$ sont des formes linéaires définies sur le fibré vectoriel $J^0_{\ell+k}TP~\times_P~J_kE$ et il est évident que $W\in\textbf{R}^0_{\ell+k}$ si et seulement si  $a^{(\lambda,\alpha)}(W)=0$ pour tout $(\lambda,\alpha)~$. La base locale $\{\epsilon^{\beta}_i\}$ de $J^0_{\ell+k}TP|U$ se remonte, par image réciproque, en une base locale du fibré $J^0_{\ell+k}TP~\times_P~\alpha^{-1}_k(U)~$. Écrivons

\begin{equation*}
a^{(\lambda,\alpha)}=\sum~a^{(\lambda,\alpha)}_{(\beta,i)}(\epsilon^{\beta}_i)^*
\end{equation*}

\vspace{2 mm}

\noindent
où $\{(\epsilon^{\beta}_i)^*\}$ est la base duale de $\{\epsilon^{\beta}_i\}$. La matrice de $\overline{\lambda}_k$ \footnote{not to be confounded with the indices $\lambda$} par rapport aux bases locales $\{\epsilon^{\beta}_i\}$ de $J^0_{\ell+k}TP~\times_P~\beta^{-1}_k(\mathcal{U})$ et $\{\partial /\partial y^{\lambda}_{\alpha}\}$ de $V\beta^{-1}_k(\mathcal{U})\subset VJ_kE$ n'est autre que $a^{(\lambda,\alpha)}_{(\beta,i)}~$, matrice dont les coefficients sont des fonctions définies sur l'ouvert $\beta^{-1}_k(\mathcal{U})$. En plus, la dimension de $(\textbf{R}^0_{\ell+k})_Z$ est égale au co-rang de cette matrice car l'isotropie ci-dessus est le noyau simultané des formes $a^{(\lambda,\alpha)}~$. Prenons maintenant un point $W\in \beta^{-1}_{k+1}(\mathcal{U})\subset J_{k+1}E$ se projetant en \textit{Z} et indiquons par $\left(b^{(\lambda,\alpha)}_{(\beta,i)}\right)~,~\alpha\leq |\alpha|\leq k+1~,~1\leq |\beta|\leq k+1~$, la matrice de $\overline{\lambda}_{k+1}~$. Le diagramme commutatif
\begin{equation*}
J^0_{\ell+k+1}TP~\times_P~J_{k+1}E\xrightarrow{\overline{\lambda}_{k+1}} VJ_{k+1}E
\end{equation*}
\begin{equation*}
\hspace{21 mm}\downarrow\hspace{34 mm}\downarrow T\rho
\end{equation*}
\begin{equation*}
\hspace{3 mm}J^0_{\ell+k}TP~\times_P~J_kE\hspace{5 mm}\xrightarrow{\overline{\lambda}_k}\hspace{4 mm}VJ_kE
\vspace{2 mm}
\end{equation*}
montre que $T\rho\circ\overline{\lambda}_{k+1}$ se factorise à $J^0_{\ell+k}TP~\times_P~J_kE$ en le morphisme $\overline{\lambda}_k$ et par conséquent la matrice de $\overline{\lambda}_{k+1}$ se transcrit par
\[\begin{pmatrix}
  B & C\\
  A & 0
 \end{pmatrix}\]
où $A$ est la matrice de $\overline{\lambda}_k$ et les blocs $(B~C)$ sont formés par les composantes des formes $b^{(\lambda,\alpha)}~,~|\alpha|=k+1~$, le bloc $C$ étant formé par les $b^{(\lambda,\alpha)}_{(\beta,j)}$ avec $|\beta|=\ell+k+1~$. Les coordonnées étant fixées, nous pouvons identifier $(J^0_{\ell+k}TP)_z$ à un facteur direct de $(J^0_{\ell+k+1}TP)_z~$, les formes $a^{(\lambda,\alpha)}_Z$ devenant des formes linéaires sur le deuxième espace. Dans ces conditions, la surjectivité de $(\textbf{R}^0_{\ell+k+1})_W\longrightarrow (\textbf{R}^0_{\ell+k})_Z$ veut dire tout simplement que le sous-espace quotient $[b^{(\lambda,\alpha)}_W]/[a^{(\lambda,\alpha)}_Z]$ se projette injectivement dans l'espace $[\epsilon^{\beta}_i (z)^*]_{|\beta|=k+1}$ où $[~~]$ indique le sous-espace engendré. Or, le rang de cette projection est donné par le rang du bloc $C_W$ et par conséquent, la surjectivité se traduit par la condition suivante: Le rang, au point $W$, de la matrice ci-dessus est égal au rang de
\[\begin{pmatrix}
  0 & C\\
  A & 0
 \end{pmatrix}\]
Puisque le rang de $C$ peut varier, cette condition ne peut pas être traduite localement par l'annulation d'un certain nombre de déterminants. Remarquons pourtant que $(\mathfrak{g}_{\ell+k+1})_W$ est le noyau simultané des restrictions, au symbole total, des formes $b^{(\lambda,\alpha)}~,~|\alpha|=k+1~$. Il en résulte que le rang de $C_W~$, égal à la codimension de $(\mathfrak{g}_{\ell+k+1})_W~$, est constant au dessus d'une trajectoire $\textbf{F}_0$ de $E$ et par conséquent, on peut, dans l'image réciproque $\beta^{-1}_{k+1}(\textbf{F}_0)~$, déterminer localement l'égalité des rangs de ces deux matrices par la nullité d'un certain nombre de déterminants. L'équation $\Theta_{k+1}$ est donc fermée dans $\beta^{-1}_{k+1}(\textbf{F}_0)~$. En particulier, $\Theta_{k+1}(\textbf{F}_0)$ est fermée dans $J_{k+1}\textbf{F}_0~$. 

\vspace{4 mm}

\noindent
Les raisonnements ci-dessus reviennent de fait à un long procéssus d'élimination. En effet, si l'on écrit explicitement les équations linéaires, portant sur les composantes $\xi^i_{\beta}$ d'un jet $j_{\ell+k+1}\xi(z)\in J^0_{\ell+k+1}TP~$, qui définissent le sous-espace $(\textbf{R}^0_{\ell+k+1})_W~$, alors la projection de cette isotro-pie dans $(\textbf{R}^0_{\ell+k})_Z$ aura, outre les équations définissant $(\textbf{R}^0_{\ell+k})_Z~$, toutes celles qui peuvent s'obtenir à partir des équtions $b^{(\lambda,\alpha)}=0~,~|\alpha|=k+1~$, par élimination d'inconnues $\xi^i_{\beta}~,~|\beta|=\ell+k+1~$. Or, cette élimination consiste à prendre les relations linéaires et linéairement indépendantes entre les lignes de $C$ c'est-à-dire les combinaisons linéaires qui annulent les lignes de $C$ et retranscrire ces mêmes relations à l'aide des lignes de $B$. On obtiendra ainsi les équations supplémentaires qui définissent la projection de $(\textbf{R}^0_{\ell+k+1})_W~$. La condition d'égalité des rangs de ces deux matrices veut dire précisément que ces équations supplémentaires sont des combinaisons linéaires des équations définissant $(\textbf{R}^0_{\ell+k})_Z~$, c'est-à-dire les équations dont les coefficients sont donnés par les lignes de $A$. On retrouve ainsi (vraisemblablement) les arguments de \cite{Pomme1973}, p.297, qui sont faits à l'ordre 1 et pour les presque-structures modelées sur un pseudogroupe de Lie transitif.  

\vspace{4 mm}

\newtheorem{close}[PropositionCounter]{Proposition}
\begin{close}
When $(\textbf{F}_0,\pi,P,p)$ is the prolongation space subordinate to an orbit $\textbf{F}_0$ contained in E then the equation $\Theta_{k+1}(\textbf{F}_0)$ in $J_{k+1}\textbf{F}_0$ is a closed sub-set locally defined by the vanishing of a certain number of determinants namely, those that establish the equality of ranks of the previously indicated matrices. Moreover, the rank of the first matrix is equal to the fibre co-dimension of $\textbf{R}^0_{\ell+k+1}|J_{k+1}\textbf{F}_0$ and the rank of the sub-matrix C is equal to the co-dimension of $\mathfrak{g}_{\ell+k+1}(\textbf{F}_0)$.
\end{close}

\vspace{4 mm}

\noindent
A prolongation space $(E,\pi,P,p)$ of order $\ell$ is said to be \textit{homogeneous} when $\Gamma_{\ell}$ or, equivalently, $\Pi_{\ell}P$ operate transitively on $E$. Examples of such spaces are the total spaces of $\ell-th$ order $\textbf{G}-$structures as well as those of almost-structures modeled on a transitive Lie pseudo-group and all the previous results transcribe replacing $\textbf{F}_0$ by $E~(=\textbf{F}_0)$. Moreover, the only proper subordinate prolongation spaces are the infinitesimal orbits in $E$. We terminate this section by resuming all the data concernig formal transitivity.

\vspace{4 mm}

\newtheorem{open}[TheoremCounter]{Theorem}
\begin{open}
Let $(E,\pi,P,p)$ be a finite prolongation space of order $\ell$ and $(\textbf{F}_0,\pi,P,p)$ the prolongation space subordinate to an orbit $\textbf{F}_0$ contained in $E$. Let us denote by $\eta_2=\eta_2(\textbf{F}_0)$ the integer where after the symbol $\mathfrak{g}_{\ell+k}(\textbf{F}_0)~,~k\geq\eta_2$ becomes $2-$acyclic and by $\eta_{\infty}$ the integer where after this symbol becomes involutive. Under these conditions:

\vspace{2 mm}

i) The non-linear equations $\mathcal{R}_{\ell+k+h}(S)~,~h\geq 0~$, associated to every solution S of an $\textbf{F}_0-$admissible fundamental equation of order $k+1\geq\eta_2+1~$, are transitive (in $(\alpha(S)$) Lie sub-groupoids that further are closed, locally trivial and regularly embedded in $\Pi_{\ell+k+h}\alpha(S)$, $\mathcal{R}_{\ell+k+h+1}(S)=p\mathcal{R}_{\ell+k+h}(S)$ and the morphism

\begin{equation*}
\mathcal{R}_{\ell+k+h+1}(S)\longrightarrow\mathcal{R}_{\ell+k+h}(S)
\end{equation*}

\vspace{2 mm}

\noindent
is a submersion. If, moreover, the isotropy groups $(\mathcal{R}^0_{\ell+k+h}S)_y$ at a point $y\in\alpha(S)$ are all connected (or else project one upon the other) then the above submersions become surjective and the equation $\mathcal{R}_{\ell+k}(S)$ will be formally integrable. This being the case, the $k-th$ order prolongation of the structure S is also formally transitive above the open set $\alpha(S)$.

\vspace{2 mm}

ii) The $\alpha-$connected components $\mathcal{R}_{\ell+k+h}(S)_0$ also verify the same regularity properties, the transitivity taking place in each connected component of $\alpha(S)$ hence everywhere when S is connected. Finally, the prolongation-projections properties stated above also transcribe for the connected components namely,  $\mathcal{R}_{\ell+k+h+1}(S)_0=p\mathcal{R}_{\ell+k+h}(S)_0$ and

\begin{equation*}
\mathcal{R}_{\ell+k+h+1}(S)_0\longrightarrow\mathcal{R}_{\ell+k+h}(S)_0
\end{equation*}

\vspace{2 mm}

\noindent
is a submersion or, in other terms, the equation $\mathcal{R}_{\ell+k}(S)_0$ is formally integrable and becomes involutive when $k\geq\eta_{\infty}~$.
\end{open}

\vspace{4 mm}

\noindent
The proof is essentially provided by all the preceding \textit{sorites}.

\vspace{2 mm}

\noindent
Let us next observe that the fibres of the fibration

\begin{equation*}
\mathcal{R}^0_{\ell+k+h+1}(S)\longrightarrow\mathcal{R}^0_{\ell+k+h}(S) 
\end{equation*}

\vspace{2 mm}

\noindent
are affine sub-spaces. Consequently, if for some integer \textit{h} the isotropy $\mathcal{R}^0_{\ell+k+h}(S)_y$ is connected, the same will occur for all higher order isotropies at that point and the transitivity of $\mathcal{R}_{\ell+k+h}(S)$ will further imply that $\mathcal{R}^0_{\ell+k+h}(S)_{y'}$ is also connected at any other point $y'\in\alpha(S)~$.  

\vspace{4 mm}

\newtheorem{closed}[CorollaryCounter]{Corollary}
\begin{closed}
The given data and the hypotheses being those of the theorem, let us assume further that $\mathcal{R}^0_{\ell+k}(S)_y$ be connected at a point $y\in\alpha(S)$. Then $\mathcal{R}_{\ell+k}(S)$ is formally integrable and the $k-th$ order prolongation of any solution S of $\textbf{F}_{k+1}$ is a formally transitive structure.  
\end{closed}

\vspace{4 mm}

\noindent
We thus see that the infinitesimal formal transitivity properties imply the same properties for the finite formal transitivity as long as we restrict the finite order equivalences to the $\alpha-$connected components of the non-linear Lie equations associated to the structure \textit{S}. Whenever one of the isotropy groups is connected, all the other goups are also connected and the formal transitivity is assured without any restriction.

\vspace{4 mm}

\noindent
We can equally define the $(k+1)-st$ order equations $\tilde{\Theta}_{k+1}$ by considering the set of all $Z\in J_{k+1}E$ such that the morphism 

\begin{equation*}
(\mathcal{R}^0_{\ell+k+1})_Z~\longrightarrow~(\mathcal{R}^0_{\ell+k})_Z
\end{equation*}

\vspace{2 mm}

\noindent
is surjective. These equations are obviously invariant by the action of $\Gamma_{\ell+k+1}$ and $\Pi_{\ell+k+1}P$ and $\tilde{\Theta}_{k+1}\subset\Theta_{k+1}~$. Moreover, every solution of a finite fundamental equation of order $k+1$  contained in $\tilde{\Theta}_{k+1}$ verifies the following two properties: 

\vspace{4 mm}

a) $\mathcal{R}_{\ell+k+1}(S)\longrightarrow\mathcal{R}_{\ell+k}(S)$ is surjective and

\vspace{2 mm}

b) $\mathcal{R}_{\ell+k}(S)$ is homogeneous of order \textit{k}.

\vspace{4 mm}

\noindent
The fact that \textit{S} is a solution of a fundamental equation implies (Lemma 8.) that both equations are "good" Lie sub-groupoids. Taking a finite trajectory $\textbf{F}_0$ in \textit{E} as well as the corresponding subordinate prolongation space, we shall obtain all the structures admitting a formally transitive prolongation - at least each connected component of such structures - as the solutions of $\textbf{F}_0-$admissible fundamental equations (\textit{i.e.}, contained in  $\tilde{\Theta}_{k+1}\cap J_{k+1}\textbf{F}_0~$, in applying the non-linear version of Quillen's Theorem at the level of $2-$acyclicity. Finally, if we ascend to the level of involutivity, we shall obtain involutive equations and consequently transitive structures, for short, in the analytic case. Besides, it is clear that the solutions of these fundamental equations are, up to a prolongation, infinitesimally formally transitive and that the preceding results concerning the properties of $\textbf{F}_0-$admissible fundamental equations can be entirely transcribed when the admissibility is taken with respect to $\tilde{\Theta}_{k+1}$ since the conditions (a) and (b) for the non-linear equations imply the corresponding conditions for the linear equations. 

\vspace{4 mm}

\noindent
Unfortunately, The equations $\tilde{\Theta}_{k+1}$ are not susceptible of a good analytical description since the surjectivity of the isotropy groups is more likely a global topological problem involving connected components once this surjectivity is verified at the Lie algebra level \textit{i.e.}, when the linear sequence
\begin{equation*}
(\textbf{R}^0_{\ell+k+1})_Z\longrightarrow (\textbf{R}^0_{\ell+k})_{\rho_k Z}\longrightarrow 0
\end{equation*}
is exact.

\noindent
We finally remark that all the preceding definitions and results are based on the general pseudo-group $\Gamma(P)$ of \textit{P} and the groupoids $\Pi_{\ell+k}P$. We can however and without any additional difficulties by simply adapting the definitions, re-write the entire section in the restricted context though it seems unavoidable, due to the nature of the various concepts, to consider no other but the transitive pseudo-groups. In fact, when $\Gamma$ is a transitive Lie pseudo-group of order $k_0$ operating on the manifold \textit{P} and $\mathcal{L}$ the corresponding infinitesimal pseudo-algebra, we can replace the general equivalence problem by the restricted one, the fundamental equations then becoming the finite or infinitesimal trajectories of the prolonged pseudo-groups or pseudo-algebras $\Gamma_{\ell+k}$ and $\mathcal{L}_{\ell+k}~$. In much the same way, we shall replace, in the notions of homogeneity and transitivity, the general equivalence by the restricted one limiting ourselves to the elements of $J_{\ell+k}\Gamma$ and $J_{\ell+k}\mathcal{L}~$.

\section{Formally transitive structures}

Let $(E,\pi,P,p)$ be a finite prolongation space of order $\ell$ that we assume, from now on, homogeneous or else we restrict our attention to a subordinate prolongation space relative to an orbit of order zero.

\vspace{4 mm}

\noindent
We next consider a formally homogeneous structure \textit{S} of species \textit{E} and will say that another structure \textit{S'} is formally equivalent to \textit{S} when, for all integers $k\geq 0$ and for any pair of points $y\in\alpha(S)$ and $y'\in\alpha(S')$, there exists $Y\in\Pi_{\ell+k}$ such that $Y(j_k S(y))=j_k S'(y')$ or, in other terms, when the two jets belong to the same $k-th$ order finite orbit namely, the orbit containing $im~j_k S~$. The jet \textit{Y} will be named a $k-th$ order finite equivalence. Further, we shall say that the formal equivalence is \textit{dominated} when there exists an integer $\nu_0$ such that every finite equivalence \textit{Y} of order $k\geq\nu_0$ is "dominated" by a finite equivalence $\tilde{Y}$ of order $k+1$ where $\rho_k\tilde{Y}=Y$. This being so, every finite equivalence of order $k\geq\nu_0$ is induced, at order \textit{k}, by an infinite order equivalence namely, by an infinite jet $(Y_h)_{h\geq 0}$ where $(Y_h)$ is a finite equivalence of order \textit{h} and $\rho_h Y_{h+1}=Y_h~$. Such an infinite jet is, in view of \textit{the Theorem of Borel}, the infinite jet $j_{\infty}\varphi(y)$ of a local diffeomorphism $\varphi$ of the manifold \textit{P}. The notion of dominated formal equivalence is actually the notion of equivalence most often adopted in the study of $\textbf{G}-$structures. The recurrent construction of principal bundles whose elements are the $k-th$ order equivalences of the flat model with the $\textbf{G}-$structure, this being achieved by requiring the nullity of the successive structure tensors or else by equivalent conditions imposed on the fundamental forms where after culminating in the desired dominated formal equivalence 

\vspace{4 mm}

\noindent
We shall say that two germs of structures of species \textit{E} are formally equivalent when they admit formally equivalent representatives, the dominated equivalence being defined similarly with the help of representatives. The study of the formal equivalence of formally infinitesimally transitive structures is a consequence of the Corollary 15.

\newtheorem{formal}[TheoremCounter]{Theorem (of formal equivalence)}
\begin{formal}
Let S be a connected infinitesimally formally transitive structure of type E and $\eta_2=\eta_2(S)$ the integer where after the symbol $\mathfrak{g}_{\ell+k}(S)~,~k\geq\eta_2~$, is $2-$acyclic. In order that a structure S' be formally equivalent to S it is necessary and sufficient that S' be a solution of the fundamental equation $\textbf{F}_{\eta_2+1}$ that contains $im~ j_{\eta_2+1}S~$. Such a structure S' is infinitesimally formally transitive and its symbol $\mathfrak{g}_{\ell+k}(S')$ has, for all $k\geq 0~$, the same homological properties of the initial one, the formal equivalence being a consequence of the finite equivalence of order $\eta_2+1$. Moreover, if S or one of its prolongations is formally transitive, that will notably take place when the non-linear isotropy groups $(\mathcal{R}^0_{\ell+k}S)_y$ at a point $y_0$ project, for $k\geq k_0~$, one upon the other or else when one of these groups is connected, then the same will hold for S' and the formal equivalence will be of dominated type.
\end{formal}

\noindent The structure \textit{S} being connected, it is also homogeneous at any order and consequently $im~j_{\eta_2+1}S$ is contained in a fundamental equation $\textbf{F}_{\eta_2+1}$ that is admissible due to the $(\eta_2+1)-st$ order infinitesimal transitivity of \textit{S}. The uniformity of the homological properties of the symbols implies that  $\mathfrak{g}_{\ell+k}=\mathfrak{g}_{\ell+k}(S)$ is $2-$acyclic for $k\geq\eta_2$ and the Corollary 15. then shows that \textit{S} and \textit{S'} are, for all $k\geq 0~$, solutions of $p_k\textbf{F}_{\eta_2+1}$ which is the admissible fundamental equation, orbit of $im~j_{\eta_2+k+1}S~$. The property relative to the symbols $\mathfrak{g}_{\ell+k}(S')$ is an immediate consequence of the Lemma 5. Finally, the part concerning the formal transitivity of \textit{S'} as well as the dominated equivalence result by an easy argument on the isotropy groups. We shall remark, for that matter, that the set of all the $k-th$ order equivalences with a given source \textit{y} and a given target \textit{y'} is a homogeneous space of the isotropy groups $(\mathcal{R}^0_{\ell+k}S)_y$ and $(\mathcal{R}^0_{\ell+k}S)_{y'}$.

\vspace{4 mm}

\noindent
$\textbf{Remark.}$ The connectedness hypothesis of the structure \textit{S} can be replaced by the formal homogeneity or even by the homogeneiy of order $\eta_2+1~$, in such a way as to assure that $im~j_{\eta_2+1}S$ be contained within an orbit.

\newtheorem{foolish}[CorollaryCounter]{Corollary}
\begin{foolish}
Let S be a connected structure or a formally homogeneous one whose $\mu-th$ prolongation is infinitesimally transitive an let $\eta_2$ be the integer from where on all the symbols $\mathfrak{g}_{\ell+k}(S)~$, $k\geq\eta_2~$, become $2-$acyclic. In order that S' be formally equivalent to S, it is necessary and sufficient that S' be a solution of the fundamental equation $\textbf{F}_{k+1}~$, $k=max~\{\mu,\eta_2\}~$, that contains $im~j_{k+1}S~$. The $\mu-th$ prolongation of such a structure S' is infinitesimally formally transitive. If, further, one of the prolongations of S is formally transitive, the same will occur with S' and the formal equivalence will be dominated.  
\end{foolish}

\vspace{2 mm}

\noindent
The formal problem being resolved, we can proceed with our inquiry and look at what really matters namely, \textit{the local equivalence} of the two structures. But here there is no way out other than to integrate the fundamental equations.

\vspace{4 mm}

\noindent
As for non-transitive structures, we should point out that much can be said in the case where the pseudo-group of local automorphisms or the pseudo-algebra of infinitesimal automorphisms have regular orbits and, more precisely, these orbits being distributed within a regular foliation. We can even assume that their prolongations to \textit{E} provide as well regular foliations. In this case we can define, at least locally, quotient prolongation spaces and quotient structures and, addapting accordingly the definitions, we shall be able to rewrite most of what was stated in this and the previous sections and consequently study the local as well as the global equivalence problems. We should also mention that Pradine's \textit{holonomy groupoid} is certainly the main tool that will hopefully determine the equivalences.

\vspace{4 mm}

\noindent
We finally say a few words on a rather outstanding example of a transitive and homogeneous structure, one of the many $Chef-d'O\!euvres$ given to us by Élie Cartan (\cite{Cartan1910}, \cite{Cartan1914}). We are referring here to the \textit{Systèmes de Pfaff en Drapeau} or, in celtic, \textit{Flag Systems}.

\vspace{4 mm}

\noindent
These are very special non-integrable Pfaffian systems enjoying the property of "slowly increasing their manifestation" of non-integrability or, in other terms, the successive derived systems decrease their dimensions just by one unit and, of course, terminate by the null system. The most remarkable about these systems is the fact that they provide a vivid example of one of Cartan's \textit{merihedric (mériédrique)} prolongation and equivalence processes. Since all this can already be learnt in \cite{Kumpera1982}, \cite{Kumpera2002} and \cite{Kumpera2013}, we shall only tell "a nice story" about them. First of all, and this is the main fact, the automorphism groups and algebras of arbitrary Flag Systems are all canonically isomorphic and, more specifically, they are all canonically isomorphic to the \textit{Darboux automorphisms} emanating from the equation $\bf{dx^2-x^3dx^1=0}~$. We can however study and classify these structures with the help of the isotropy groups and algebras of arbitrary order. These are intrinsically defined objects and the classification comes out upon looking at the differences in the \textit{co-ranks} of these (Fréchet infinite dimensional) groups and algebras at least for those flags that the author in \cite{Kumpera2014} calls \textit{elementary}. When the systems are non-elementary, their classification becomes much more delicate and a new phenomenon does occur namely, the appearance of continuous or even differentiable deformations (variations) of a given flag via non-equivalent flags. However, we can also cope with such phenomena by considering again the isotropy algebras as discrete (finite number of variables) moduli. Most authors still employ the rather clumsy and outdated \textit{pseudo-normal forms} introduced, way back in 1980, by the present author though the above terminology is due to others. At present, we can do no better than apologize for having considered these silly and uncouth objects and where several authors come forth with fathom deep proofs that nobody will ever understand except, and with all the blessing and mercy of the All Mighty, those authors that wrote them.

\vspace{4 mm}

\noindent
And what about the totally intransitive or rigid structures. For this, we shall be forced to "\textit{admirer le beau et fragile papillon qui ouvre ses ailes à Aix-en-Provence et provoque non le dépliement des ailes d'un autre papillon mais en fait une tsuname à Fukushima}". This being a rather long subject we shall not even try to say a word about it hoping that other authors give us much to read and enjoy.

\bibliographystyle{plain}
\bibliography{references}
\end{document}